\documentstyle[amsfonts,amssymb]{article}

\newtheorem{theorem}{Theorem}[section]
\newtheorem{lemma}[theorem]{Lemma}

\newtheorem{corollary}[theorem]{Corollary}

\newtheorem{proposition}[theorem]{Proposition}
\newtheorem{example}[theorem]{Example}

\title{\bf Global solutions for semilinear Klein-Gordon equations in    
FLRW spacetimes} 

\author{{\bf Anahit Galstian${}^{\small \mbox{\rm \scriptsize 1}}$  and Karen Yagdjian${}^{\scriptsize \mbox{\rm \scriptsize 2}}$ } }

\begin{document}

\date{}

\maketitle

\thispagestyle{empty}

\begin{center}
\noindent
{
{\small \textsl{Department of Mathematics,
University of Texas-Pan American,  \\ 1201 W.~University Drive,  Edinburg, TX 78539, USA\\
{ ${}^1$agalstyan@utpa.edu }}\,
 ${}^2$\textsl{ yagdjian@utpa.edu}}}
\end{center}

\begin{abstract}
We consider   waves, which obey the semilinear Klein-Gordon equation,  propagating in the Friedmann-Lema$\rm \hat{i} $tre-Robertson-Walker spacetimes.
The equations   in the de~Sitter and  Einstein-de~Sitter spacetimes  are the important particular cases.
 We show the global in time existence   in the energy class of solutions of the Cauchy problem.

\medskip

\noindent
{\bf Keywords\,\, } Klein-Gordon equation $ \cdot $   Einstein-de~Sitter universe   $ \cdot $  de~Sitter  universe  $ \cdot $ global solutions
\end{abstract}

\section{Introduction}
\setcounter{equation}{0}
\renewcommand{\theequation}{\thesection.\arabic{equation}}

In this article we  consider the   Klein-Gordon equation in the spacetimes belonging to some family of  the  Friedmann-Lema$\rm \hat{i} $tre-Robertson-Walker spacetimes (FLRW spacetimes).
In the FLRW  spacetime, one can choose coordinates so that the metric has the form $
ds^2=-dt^2+a^2(t)d \sigma ^2$. (See, e.g., \cite{Hawking}.) This family includes, as a particular case,  the metric
\begin{equation}
\label{metric_l}
ds^2= -  \, dt^2+  t^\ell  \sum_{i,j=1,\ldots,n}\delta _{ij}  dx^idx^j\,,
\end{equation}
where $\delta _{ij} $ is the Kronecker symbol and $\ell= \frac{4}{n\gamma}$.  The function $a (t)$ is the scaling factor. The time dependence of the function $a(t)$ is  determined by  the Einstein's field equations for gravity, which  for the perfect fluid  imply  
\begin{eqnarray*}
&  &
\dot \mu  =-3(\mu  +p)\frac{\dot a}{a}\,, \qquad
\frac{\ddot  a}{a} =-\frac{4\pi }{3} (\mu  +3p)\,, \qquad 
\label{Feq}
\left( \frac{\dot  a}{a} \right)^2 = \frac{8\pi }{3}  \mu   - \frac{K}{a^2}\,,
\end{eqnarray*}
where $\mu  $ is the proper energy density,   $p $ is pressure, and  $K  $ is  spatial curvature. 
The last equations give a   differential equation for $a=a(t)$ if an equation of state 
(equation for the pressure) $p=p(\mu  )$
is known.  
For pressureless,  $p = 0$, matter distribution in the universe and vanishing
spatial curvature, $K = 0$,  the
solution to that  
equation is
\begin{eqnarray*}
&  &
a(t)= a_0 t^{{2}/{3}}\,,
\end{eqnarray*}
where $ a_0>0$   is a  constant.   The universe expands, and its expansion decelerates since 
$\ddot{a}   < 0$.  In the radiation dominated universe, the equation of state is $p=\mu /3 $ and, consequently, $a(t)= a_0 t^{{1}/{2}} $.  The equation of state $p=(\gamma -1)\mu  $, which includes those two cases of the matter-  and   radiation- dominated universe,  implies that in order to have a non-negative pressure for a positive density, 
it must   be assumed that 
 $\gamma \geq 1 $ for the physical space with $n=3$ (see \cite{Choquet-Bruhat_book} p.122).  The spacetime with   $\gamma =1 $
and $n=3$  is called the Einstein-de~Sitter universe.
We reveal in this paper the significance  of the   restriction $\gamma \geq 1 $ on the range of $\ell$; in fact, we show that it is  closely related to the non-growth of the  energy   and to the existence of the global in time solution of the Cauchy problem for the   Klein-Gordon equation. Another important spacetime, the so-called de~Sitter spacetime, is also a member of that family and it will be discussed as well.
\medskip

In  quantum field theory    the matter fields are described by a function  $\psi $  that must satisfy  equations of motion.
In the case of a massive scalar field, the equation of motion is   the  semilinear Klein-Gordon equation generated by the metric $g$:
\begin{eqnarray}
\label{CKGE}
\frac{1}{\sqrt{|g(x)|}}\frac{\partial }{\partial x^i}\left( \sqrt{|g(x)|} g^{ik} (x)\frac{\partial \psi   }{\partial x^k} \right) = m^2 \psi   +   V'_\psi (x,\psi  ) \,.
\end{eqnarray}
In physical
terms this equation describes a local self-interaction for a scalar particle.  A typical example of a potential function  would be  $V(\phi )=\phi ^4$.
The semilinear equations are also commonly used models for general nonlinear problems. 
\smallskip

To motivate our approach, we first consider the  covariant Klein-Gordon equation  in the metric (\ref{metric_l}), which can be written in the global coordinates  as follows
 \begin{eqnarray*}
    \psi _{tt}
-  t^{-\ell } \Delta     \psi  + \frac{ n \ell}{ 2 t }       \psi_t + m^2 \psi
+ V'_\psi (x,t,\psi  )& = &
 0 \,.
\end{eqnarray*}
Here $x \in {\mathbb R}^n$, $t \in {\mathbb R}$, and $\Delta  $ is
 the Laplace operator on the flat metric, $\Delta  := \sum_{j=1}^n \frac{\partial^2 }{\partial x_j ^2} $.

We study the Cauchy problem with the data prescribed at some positive time $t_0  $
\begin{eqnarray}
\label{IC}
   \psi (t_0 ,x)= \psi _0, \quad    \psi_t (t_0  ,x)= \psi _1\,,
\end{eqnarray}
and we look for the solution defined for all values of $t \in [t_0  ,\infty)$   and $x \in {\mathbb R}^n $.
We change the unknown function $\psi =t^{-\frac{ln}{4}}u $, then  for the new function $u=u(t,x) $ we obtain the equation
\begin{eqnarray*}
  u_{tt} -  t^{-\ell } \Delta    u +  M^2(t)   u + t^{n\ell /4}V'_\psi (x,t,t^{-\frac{\ell n}{4}}u )
& = &
0
\end{eqnarray*}
with the ``effective'' (or ``curved mass'')
\begin{equation}
\label{MEdS}
M^2_{EdS}(t) :=   m^2- \frac{ n \ell (  n \ell-4)}{16 t^2 }  \,.
\end{equation}
It is easily seen that  for the range $(0,\frac{4}{n}] $ of the parameter $\ell $ the curved mass is positive while its derivative is non-positive.
This is crucial for the non-increasing property of the energy and in the derivation of the energy estimate.
\medskip

Let $(V,g) $ be smooth pseudo  Riemannian  manifold of dimension $n+1 $ and $V= {\mathbb R} \times S$ with $S $ an $n$-dimensional orientable
smooth 
manifold, and $g$ be a FLRW  metric.
We restrict our attention to the case of  $n\geq 3 $ and to
the spacetime with the line element
\begin{eqnarray*}
 &  &
 ds^2= -  \, dt^2+   a^{2}(t)  \sigma \,.
 \end{eqnarray*}
Then we consider an {\it expanding universe} that means that $ \dot{a}(t)>0$.
For the metric with $ \dot{a}(t)>0$ we define the norm
\begin{eqnarray}
\label{normX}
 \parallel \psi   \parallel _{X(t)}
& :=  &
\| \psi  _{t}   \| _{L^\infty([t_0  ,t]; L^2(S))} +    \| a^{-1} (\cdot )\nabla_\sigma   \psi     \| _{L^\infty([t_0  ,t]; L^2( S))}  \\
&  &
 + \| M(\cdot )\psi   \| _{L^\infty([t_0  ,t]; L^2(S))}+   \|     \sqrt{ \dot{a}  a^{-3}  }\, \nabla_\sigma  \psi   \|_{L^2([t_0  ,t]\times S )}
\nonumber\, ,
\end{eqnarray}
where $0< t_0  <t \leq \infty $ and $M(t) \geq 0$ is a curved mass defined by:
\begin{eqnarray}
\label{M}
M^2(t)  = m^2 +\left( \frac{n}{2}-\frac{ n^2 }{ 4 } \right)\left( \frac{\dot{a}(t)}{a (t)}\right)^2  -\frac{n}{2} \frac{\ddot{a}(t)}{a (t)} \,.
\end{eqnarray}
Hence, in the classification suggested in \cite{Yagdjian-Galstian}, mass $m$ is large if the metric $g$ is  a de~Sitter metric $-  \, dt^2+   e^{2t}  dx^2 $, $x \in {\mathbb R}^n$.
Here and henceforth  $ \dot{a} (t)$ denotes  the derivative with respect to time, while   the spatial variable will be denoted $s$
in a general manifold $S$ and $x$ when $S = {\mathbb R}^n$.
 In order to describe admissible nonlinearities  we make the following definition.
\smallskip

\noindent
 {\bf Condition ($\mathcal L$)}. {\it The function $F=F(s,u)$,\,\, $F\,:\, S\times {\mathbb R}\longrightarrow {\mathbb R}$ is said to be Lipschitz continuous  in $u$ with exponent $\alpha  $, if  there exist $\alpha \geq 0$ and $C>0$ such that}
\[
 \hspace{-0.4cm} |  F(s,u )- F(s,v )|  \leq C
|  u -  v  |
\left( | u  |^\alpha
+ |  v |^\alpha \right) \quad \mbox{\sl for all} \,\,  u,v \in  {\mathbb R}, \,\,  x \in  S\,.
\]
For the continuous function $\Gamma \in C([t_0, \infty)) $ denote by $C_{a,\Gamma ,\alpha _0}(T)$ and  $C_{a,\Gamma ,\alpha _0}^{(-1)}(r)$    the function
\begin{eqnarray*}
&    &
C_{a,\Gamma ,\alpha _0}(T):= \left( \int_{t_0}  ^T  \left( \frac{a(t)}{ \dot{a}(t)} \right) ^{\frac{n\alpha _0  }{4-n\alpha _0} }  \left| \Gamma (t)\right|^{\frac{4}{4-n\alpha _0}} \,dt \right)^{\frac{4-n\alpha _0}{4}}    ,\qquad 0 <  \alpha _0 < \frac{4}{n}  \,.
\end{eqnarray*}
and its inverse, respectively. 
The main result of this paper is the following theorem.
\begin{theorem}
\label{T1.1}
Assume that \, $n=3,4$ \,and that the metric $g $ is $
g= -  \, dt^2+   a^{2}(t)  \sigma$. 
Suppose also that  $ m>0$ and that  there is a positive number $c_0$ such that  the real-valued positive function $a=a(t) $ satisfies
\begin{eqnarray}
\label{a}
&  &
{a} (t) > 0, \quad  \dot{a} (t) > 0 \quad  \mbox{for all }   \quad t  \in [t_0  ,\infty)\,,\\
\label{Ma}
&  &
M(t)
  > c_0> 0 ,     \quad  \dot{M} (t) \leq 0 
  \quad \mbox{ for all}   \quad t \in [t_0  ,\infty)  \,\, \,.
\end{eqnarray}
Consider the Cauchy problem for the equation (\ref{CKGE}) with 
the derivative of potential function $V'_\psi (s, t, \psi )= -\Gamma (t) F(s,\psi )$ such that $ F$ is   Lipschitz continuous with exponent $\alpha  $, $F(s,0)=0$ for all $s \in S $, and  either
\begin{eqnarray}
\label{GammaB}
&    &
 | \Gamma (t) | \leq C_\Gamma  \frac{  \dot{a}(t) }{a(t)}    \quad  \mbox{for all }   \quad t  \in [t_0  ,\infty)\,,
\end{eqnarray}
where $C_\Gamma $ is a constant independent of $t$, or, there is $\alpha _0 $ such that
\begin{eqnarray}
\label{GammaU}
&    &
C_{a,\Gamma ,\alpha _0}(\infty) < \infty  ,\qquad 0 <  \alpha _0 < \frac{4}{n}  \,.
\end{eqnarray}

If \, $ \frac{4}{n} \leq  \alpha \leq \frac{2}{n-2} $, then for every $\psi _0 \in  H_{(1)}(S)$ and  $\psi _1\in  L^2(S)$, sufficiently small initial data,  $\| \psi _0\|_{H_{(1)}(S)} + \| \psi _1\|_{L^2(S)}$, the problem (\ref{CKGE}),(\ref{IC}) has a unique solution $\psi \in C([t_0  ,\infty);H_{(1)}(S) )\cap  C^1([t_0  ,\infty);$ $L^2(S))$ and its norm $\| a ^{\frac{n}{2}}\psi \|_{X(\infty)}$ is small.
\end{theorem}

Condition (\ref{Ma})   for the norm of solutions of the equation   implies that the energy of solution is non-increasing.
In the next theorem the local existence is stated with the less restrictive conditions and with the estimate for the lifespan.
\begin{theorem}
\label{T1.3}
Suppose   that  $ m>0$ and that  there is a positive number $c_0$ such that  the real-valued positive function $a=a(t) $ satisfies (\ref{a}),
(\ref{Ma}).
Consider the Cauchy problem for the equation (\ref{CKGE}) with 
the derivative of potential function $V'_\psi (s, t, \psi )= -\Gamma (t) F(s,\psi )$ such that $ F$ is   Lipschitz continuous with exponent $\alpha  $, $F(s,0)=0$ for all $s \in S $.

If  \, $ 0 \leq  \alpha \leq \frac{2}{n-2} $, then for every $\psi _0 \in  H_{(1)}(S)$ and  $\psi _1\in  L^2(S)$ 
there exists  \,$T_1>t_0$ \, such that the problem (\ref{CKGE}),(\ref{IC}) has a unique solution $\psi \in C([t_0  ,T_1);H_{(1)}(S) )\cap  C^1([t_0  ,T_1);L^2(S))$.

The  lifespan of the solution can be estimated as follows
\[
T_1- t_0  \geq C C_{a,\Gamma ,\alpha _0}^{(-1)}(\| \psi _0\|_{H_{(1)}(S)} + \| \psi _1\|_{L^2(S)})\,,
\]
where $C$ is a positive constant independent of $T_1$, $\psi _0 $ and $\psi _1 $.
\end{theorem}
\smallskip

If the  nonlinear term has an energy conservative potential function, then in the next theorem we establish the existence of the global solution for large initial data.
\begin{theorem}
\label{T1.2}
Suppose that all conditions of Theorem~\ref{T1.3} on\,  $ n$, $\alpha  $,  and $a=a(t) $,  are satisfied, and additionally,
 \begin{equation}
\label{potential}
\frac{2}{n}\frac{a (t)}{\dot{a} (t)}V_t'(t,s,a^{-n/2}(t)w)  +  2V (t,s,a^{-n/2}(t)w)   
-   a^{-n/2}(t) w V_\psi ' (s,a^{-n/2}(t)w)\leq 0
\end{equation}
for all $(t,s,w)   \in [t_0,\infty) \times S\times {\mathbb R}$.

Then for every $\psi _0 \in  H_{(1)}({\mathbb R}^n)$ and  $\psi _1\in  L^2({\mathbb R}^n)$, the problem (\ref{CKGE}),(\ref{IC}) has a unique solution $\psi \in C([t_0  ,\infty);H_{(1)}({\mathbb R}^n)) \cap  C^1([t_0  ,\infty);L^2({\mathbb R}^n))$ and  its norm $\| a^{\frac{n}{2}}\psi \|_{X(\infty)}$ is finite.
\end{theorem}
\medskip

The hyperbolic equations in the de~Sitter 
spacetime   have permanently bounded domain of  influence. Nonlinear equations with a permanently bounded domain of  influence
were studied, in particular, in \cite{Yag_2005}. In that paper  the example of   equation, which has  a blowing-up solution for  arbitrarily small data, is given. Moreover,
it was discovered in \cite{Yag_2005} that the time-oscillation of the metric, due to the parametric resonance, can cause   blowup phenomena for   wave map type nonlinearities even for the arbitrarily small data. On the other hand in the absence of   oscillations in the metric,  Choquet-Bruhat \cite{Choquet-Bruhat-ND_2000} proved for small initial data
the global existence and uniqueness of wave maps on the 
FLRW expanding universe with the metric ${\bf g}=-dt^2+R^2(t)   \sigma   $  and    a smooth Riemannian manifold $(S,\sigma )$ of dimension $n \leq 3$, which has a 
  time independent metric $\sigma  $  and   non-zero injectivity radius,
and with $R(t) $ being  a positive increasing function  such that $1/R (t)$ is integrable on $[t_0  ,\infty)$.   If the target manifold is flat,  then the wave map equation reduces to a linear system. On the other hand, in the Einstein-de~Sitter spacetime the domain of influence is not permanently bounded.
\smallskip

In the   de Sitter space, that is, in the spacetime with the line element
\begin{eqnarray}
\label{1.6}
 &  &
 ds^2= -  \, dt^2+   e^{2t}    \sum_{i,j=1,\ldots,n}\delta _{ij}  dx^idx^j\,,
 \end{eqnarray}
the second author \cite{JMAA_2012}-\cite{Helsinki_2013}    studied the Cauchy problem for
the semilinear equation
\begin{eqnarray*}
&  &
\square_g u + m^2 u = f(u),\qquad
u(x,t_0)=\varphi_0 (x) \in H_{(s)}({\mathbb R}^n) , \,\,\, u_t(x,t_0)=\varphi_1 (x) \in H_{(s)} ({\mathbb R}^n)\,,
\end{eqnarray*}
if $s >n/2$ and $f$ is   Lipschitz continuous with exponent $\alpha  $. In \cite{JMAA_2012}-\cite{Helsinki_2013} a global existence of   small data solutions of the Cauchy problem for the semilinear  Klein-Gordon equation and systems of equations  in the de~Sitter spacetime  is proved.
It was discovered that unlike   the same problem in the Minkowski spacetime,   no restriction on the order of nonlinearity is required,
provided that a physical mass of the field belongs to some set,  $m \in (0, \sqrt{n^2-1}/2]\cup [n/2,\infty) $.
It was also conjectured  that $  (  \sqrt{n^2-1}/2,n/2 ) $ {\it is a forbidden mass interval for the small data global solvability of the Cauchy problem for all}  $\alpha \in (0,\infty) $. For $n=3$ the interval $(0,\sqrt{2}) $ is  called  the Higuchi bound  in  quantum field theory \cite{Higuchi}. The proof of the global existence in \cite{JMAA_2012}-\cite{Helsinki_2013} is based on the $ L^p-L^q$ estimates.
\smallskip

Baskin~\cite{BaskinSE} discussed small data global solutions for the
scalar Klein-Gordon equation on asymptotically de Sitter spaces,
which are   compact manifolds with boundary.   More precisely, in
\cite{BaskinSE}  the following Cauchy problem is considered for the
semilinear equation
\begin{eqnarray*}
&  &
\square_g u +m^2 u = f(u),\qquad
u(x,t_0)=\varphi_0 (x) \in H_{(1)}({\mathbb R}^n) , \,\,\, u_t(x,t_0)=\varphi_1 (x) \in L^2 ({\mathbb R}^n)\,,
\end{eqnarray*}
where mass is large, $m^2> n^2/4$, $f$ is a smooth function and
satisfies conditions $ |f(u)| \leq c|u|^{\alpha  +1}$, $ |u| \cdot
|f'(u)| \sim |f(u)|$, $f(u)-f'(u)\cdot u   \leq 0$, $  \int_0^u
f(v)dv \geq 0$, and $\int_0^u f(v)dv \sim |u|^{\alpha  +2} $ for
large $|u|$. It is also assumed that $  \alpha= \frac{4}{n-1}$. In
Theorem~1.3 \cite{BaskinSE} the existence of the global solution for
small energy data is stated. (For more references on the
asymptotically de Sitter spaces, see  the bibliography in
\cite{Baskin},
 \cite{Vasy_2010}.)
\smallskip

 Hintz and Vasy \cite{Hintz-Vasy} considered semilinear wave equations
of the form
\[
(\square_g-\lambda ) u =f+q(u,du)
\]
on a manifold $M$, where $q$ is  typically
a polynomial vanishing at least quadratically at $(0, 0)$,
 in contexts such as asymptotically de Sitter and Kerr-de Sitter spaces, as well as asymptotically Minkowski
spaces.  The linear framework in \cite{Hintz-Vasy} is based
on the b-analysis, in the sense of Melrose, introduced in this context by Vasy
to describe the asymptotic behavior of solutions of linear equations. Hintz and Vasy  have shown the small data solvability of suitable semilinear
wave and Klein-Gordon equations.
\smallskip

Nakamura   \cite{Nakamura} considered the Cauchy problem for the semi-linear Klein-Gordon equations  in de~Sitter spacetime with $n \leq 4 $, that is, with the line element
(\ref{1.6}).
The nonlinear term  is of power type for $n=3,4$, or of exponential type for $n=1,2$. For the power type semilinear term with \, $ \frac{4}{n} \leq  \alpha \leq \frac{2}{n-2} $ Nakamura \cite{Nakamura} proved  the existence of    global solutions in the energy class.
\smallskip

Ringstr\"om~\cite{Ringstrom} considered the question of future global non-linear stability in the case of Einstein's equations coupled to a non-linear scalar field. The class of potential $V(\psi )$ is restricted by the condition $V(0)>0$, $V'(0)=0$ and $V''(0)>0 $. Ringstr\"om proved that for given initial data, there is a maximal globally hyperbolic development of the data which is unique
up to isometry. The case of Einstein's equations with positive cosmological constant 
was not included  unless the scalar field is  zero. 
\smallskip

Rodnianski and Speck \cite{Rodnianski-Speck} proved the nonlinear future stability of the FLRW family of solutions to the irrotational Euler-Einstein system with a positive cosmological constant. 
More precisely, they studied small perturbations of the family of FLRW cosmological background solutions to the coupled Euler-Einstein system with
a positive cosmological constant in $1 + 3$ spacetime dimensions. The background solutions model
an initially uniform quiet fluid of positive energy density evolving in a spacetime undergoing exponentially
accelerated expansion. Their analysis shows that under the equation of state $p=(\gamma -1)\mu$, $0<\gamma -1<1/3 $,  the background metric plus fluid 
solutions are globally future-stable under small irrotational perturbations of their initial data.

\medskip

The present paper is organized as follows.  In Section \ref{S2} we derive the energy estimates. Then,
 in Section \ref{S3}, we give  the estimate for the nonlinear term. The completion of the proof of Theorems~\ref{T1.1},\ref{T1.3},\ref{T1.2} is given in Section~\ref{S4}. To illustrate our results we
discuss  some examples in that section.

\medskip

\section{The linear Klein-Gordon equation}
\label{S2}
\setcounter{equation}{0}

\subsection{FLRW universe. The effect of the damping term}

The goal of the transformation that has been used in \cite{Yagdjian-Galstian} and  will be used in this subsection is twofold; on one hand, it reduces the damping term and makes possible to prove non increase property of the energy  and, on the another hand, 
it provides the nonlinear term with some weight function, which    is decreasing in time. In fact, this is a particular case of the Liouville transformation that is used to study boundedness, stability and asymptotic behavior 
of solutions of the second order differential equations.

The line element  in the theorems for the  FLRW spacetime implies
\[
g_{00}=  g^{00}= -  1  ,\,\, g_{0j}= g^{0j}= 0, \,\,  g_{ij}=a^2(t)  \delta _{ij} (x) ,\,\,  |g|=  a^{2n}(t)  |\det \delta (x)|,\,\,            g^{ij}=  a^{-2}(t)\delta ^{ij} (x)  ,
\]
$i,j=1,2,\ldots,n$, where $\delta ^{ij} (x) \delta _{jk} (x)=\delta_{k}^i $  (Kronecker symbol) and the metric $\sigma$ in the local chart is given by $\delta _{ik} (x)$.

The linear covariant Klein-Gordon equation in that background is  $\square_g \psi
 =
m^2 \psi -f $ and in the coordinates this equation can be written as follows
\begin{eqnarray*}
 \psi _{tt}
-  \frac{1}{a^2(t)\sqrt{|\det \delta ( x)| }} \sum_{i,j=1}^n  \frac{\partial  }{\partial x^i}\left(  \sqrt{|\det \delta ( x)| }  \delta ^{ij} (x)\frac{\partial  }{\partial x^j}  \psi \right) + n \frac{ \dot{a}(t)}{a (t)}     \psi_t + m^2 \psi
& = &
   f \,.
\end{eqnarray*}
In order to eliminate the damping term $  n \frac{ \dot{a}(t)}{a (t)}     \psi_t$  we introduce  the new unknown function  $\psi =b(t) u $, then
the equation for $u$ is
\begin{eqnarray}
\label{eqb}
&  &
 u_{tt}-  \frac{1}{a^2(t)\sqrt{|\det \delta (x)| }} \sum_{i,j=1}^n  \frac{\partial  }{\partial x^i}\left(  \sqrt{|\det \delta (x)| }  \delta ^{ij} (x)\frac{\partial  u }{\partial x^j}  \right)+ \left( 2\frac{ \dot{b} (t)}{b(t)} + n \frac{ \dot{a}(t)}{a (t)}  \right) u_t
\nonumber \\
&  &
  +   \left(   n \frac{ \dot{a}(t)}{a (t)}\frac{\dot{b}(t)}{b(t)} +  \frac{\ddot{b}(t)}{b(t)}\right)  u  + m^2 u
=  \frac{1}{b(t)}f\,.
\end{eqnarray}
We look for the function $ b=b(t)$ such that the following equation is fulfilled:
\begin{eqnarray*}
&  &
2\frac{ \dot{b} (t)}{b(t)} + n \frac{ \dot{a}(t)}{a (t)} =0 \,.
\end{eqnarray*}
In particular, we can choose
\begin{eqnarray}
\label{trans}
&  &
 \displaystyle  b(t) =  a^{-\frac{n}{2}}(t)  \,.
\end{eqnarray}
Consequently,  the coefficient of  $u_t $ vanished while the coefficient of the term with $u $ of the equation (\ref{eqb}) became
\begin{eqnarray*}
&  &
n \frac{ \dot{a}(t)}{a (t)}\frac{\dot{b}(t)}{b(t)} +  \frac{\ddot{b}(t)}{b(t)}
=
\left( \frac{n}{2}-\frac{ n^2 }{ 4 } \right)\left( \frac{\dot{a}(t)}{a (t)}\right)^2  -\frac{n}{2} \frac{\ddot{a}(t)}{a (t)}\,.
\end{eqnarray*}
Hence, the Klein-Gordon equation for the function $u=u(x,t)$ can be written as follows:
\begin{eqnarray*}
 &  &
u_{tt}-  \frac{1}{a^2(t)\sqrt{|\det \delta (x)| }}  \sum_{i,j=1}^n \frac{\partial  }{\partial x^i}
\left(  \sqrt{|\det \delta (x)| }  \delta ^{ij} (x)\frac{\partial  u }{\partial x^j}  \right)\\
&  &
\hspace{4cm} +\left( m^2 +\left( \frac{n}{2}-\frac{ n^2 }{ 4 } \right)\left( \frac{\dot{a}(t)}{a (t)}\right)^2  -\frac{n}{2} \frac{\ddot{a}(t)}{a (t)}  \right)u
 =
 \frac{1}{b(t)}f\,.
\end{eqnarray*}
We denote the coefficient of the last equation by
\begin{eqnarray*}
c(t) :=  m^2 +\left( \frac{n}{2}-\frac{ n^2 }{ 4 } \right)\left( \frac{\dot{a}(t)}{a (t)}\right)^2  -\frac{n}{2} \frac{\ddot{a}(t)}{a (t)}  \,.
\end{eqnarray*}
For the  FLRW spacetime  with  \, $a(t)=t^{\ell/2}$, \, $\ell  \leq \frac{4}{n} $,  and
\[
g_{00}=  g^{00}= -  1  ,\,\, g_{0j}= g^{0j}= 0, \,\,  g_{ij}=t^\ell  \delta _{ij} (x) ,\,\,  |g|=   t^{\ell  n} |\det \delta (x)|,\,\,            g^{ij}=  t^{-\ell }\delta ^{ij} (x)   \,,
\quad
\]
$i,j=1,\ldots,n$, in accordance with \cite{Choquet-Bruhat_book} p.124,  we have
\begin{eqnarray*}
c (t) = m^2- \frac{ n\ell ( n\ell -4)}{16 t^2 }> 0 \,,  \quad    \dot{c}  (t) = \frac{ n\ell (  n\ell-4)}{8 t^3 }\leq 0
\quad  \mbox{\rm  for }   \quad \quad  m > 0, \quad   \ell  \leq \frac{4}{n}\,.
\end{eqnarray*}
The last inequalities suggest the assumptions on $a(t) $, $m$, and $n$:
\begin{eqnarray*}
c (t)
& = &
m^2 +\left( \frac{n}{2}-\frac{ n^2 }{ 4 } \right)\left( \frac{\dot{a}(t)}{a (t)}\right)^2  -\frac{n}{2} \frac{\ddot{a}(t)}{a (t)}    > 0 \quad  \mbox{\rm  for all large}   \quad t \quad  \mbox{\rm  and}   \quad   m > 0 \,,  \\
   \dot{c}  (t)
& = &
\frac{d}{d t} \left( \left( \frac{n}{2}-\frac{ n^2 }{ 4 } \right)\left( \frac{\dot{a}(t)}{a (t)}\right)^2  -\frac{n}{2} \frac{\ddot{a}(t)}{a (t)}   \right)   \leq 0  \quad  \mbox{\rm  for all large}   \quad t   \,.
\end{eqnarray*}
Thus, in order to study the equation (\ref{CKGE}) we can  first consider the following linear equation
\begin{eqnarray*}
 &  &
u_{tt}-  \frac{1}{a^2(t) } \Delta_\sigma u
+M^2(t) u  =   \frac{1}{b(t)}f\,,
\end{eqnarray*}
with the curved mass $M(t) $, $M^2(t)  = c(t) $, of (\ref{M}), 
and derive in the next subsection for the solutions the energy estimates.
Here $ \Delta _\sigma $ is a Laplace-Beltrami operator in the metric $\sigma $.

\subsection{Energy estimate}

\renewcommand{\theequation}{\thesection.\arabic{equation}}

In this subsection we show that for the large physical mass $m$ the expansion property of the de Sitter metric leads, via transformation (\ref{trans}),  to the dissipative effect for the Klein-Gordon equation. Eventually   the  dissipation   prevents the blow up  of the solution of  the nonlinear equation.

First we consider the solution $u=u(t,x) $ of the equation without source term, $ f=0$,
\begin{eqnarray*}
 &  &
 u_{tt}-  \frac{1}{a^2(t)} \Delta _\sigma
+M^2(t) u  =   0\,.
\end{eqnarray*}
For a Riemannian manifold $(S,\sigma )$ we denote by $\nabla_\sigma $ the covariant gradient and by $ d\mu _\sigma $ the volume element in the metric $\sigma  $. The Sobolev space $W^p_s (S)$ is a Banach space with the norm
\[
\| u \|_{W^p_s  (S)} := \left(\int_S \sum_{0\leq |\alpha | \leq s} |\partial ^\alpha  u|^p\, d\mu _\sigma \right) ^{1/p}, \quad 1 \leq p <\infty\,.
\]

We define the energy of the solution $u=u(t,s) $  by
\begin{eqnarray}
\label{energy}
E(t)  & :=   &
\frac{1}{2} \| u_{t}\|^2_{ L^2( S )}+ \frac{1}{2}a^{-2}(t)\| \nabla_\sigma  u \|^2_{ L^2( S )}
 +  \frac{1}{2}M^2(t)   \|      u   \|^2_{ L^2(S )} \,.
\end{eqnarray}
Then 
\begin{eqnarray*}
\frac{d}{d t }  E(t )
& = &
  \frac{1}{2}  ( a^{-2}(t) )_t\| \nabla_\sigma  u \|^2_{ L^2( S )}
+  \frac{1}{2}  (M^2 (t) )_t   \|      u   \|^2_{ L^2(S )} \leq 0\,.
\end{eqnarray*}
Integration gives
\begin{eqnarray*}
E(t)
& -  &
\int_{t_0}  ^t \Bigg[  \frac{1}{2}  ( a^{-2}(\tau ) )_\tau \| \nabla_\sigma  u \|^2_{ L^2( S )}
+  \frac{1}{2}  (M^2 (\tau ) )_\tau    \|      u   \|^2_{ L^2(S )}  \Bigg] \,d \tau  = E(t_0  )\,.
\end{eqnarray*}
In particular,  due to the  assumptions   on  $a(t)$  and  $M(t)$  we obtain
$
E(t)\leq  E(t_0  )$, 
that is,
\begin{eqnarray*}
&  &
 \| u_{t}\|^2_{ L^2( S )}+ a^{-2}(t)\| \nabla_\sigma  u \|^2_{ L^2( S )}
 +  M^2(t)   \|      u   \|^2_{ L^2(S )}   \\
&  \leq &
 \| u_{t}(t_0  ) \|^2_{ L^2( S )}+  a^{-2}(t_0  )\| \nabla_\sigma  u(t_0  ) \|^2_{ L^2( S )}
 +   M^2(t_0  )   \|      u (t_0  )  \|^2_{ L^2(S )}
\,.
\end{eqnarray*}
We    also have  $0< M_0 \leq M(t) \leq M_1 < \infty$ for all $t \geq t_0$ with some constants $ M_0$, $M_1 $.

Consider now the solution $u $ of the equation
\begin{eqnarray}
\label{5.7g}
 u_{tt}
-   a^{-2}(t)   \Delta _\sigma   u      +M^2(t)  u =g
\end{eqnarray}
with the source term $g$. Here and henceforth, if $A$ and $B$ are two non-negative quantities, we use $A \lesssim  B$ to denote the statement that $A\leq CB $ for some absolute constant $C>0$.
\begin{proposition}
Assume that conditions (\ref{a}) and (\ref{Ma}) are fulfilled. Then,  the solution $ u=u(t,s)$ of the equation (\ref{5.7g}) satisfies the following  estimate
\begin{eqnarray*}
 &  &
    \| u_{t}   \| _{L^\infty([t_0  ,t]; L^2(S))} +  \| a^{-1} (t)\nabla u  (x,t ) \| _{L^\infty([t_0  ,t]; L^2(S))}
 + \| M (t)u  (x,t ) \| _{L^\infty([t_0  ,t]; L^2(S))}  \nonumber \\
&  &
+  \|   \sqrt{ \dot{a}  a^{-3}  } \,\nabla u   \|_{L^2([t_0  ,t]\times S)}   +  \| \sqrt{ |\dot{c} | } \, u   \| _{L^2([t_0  ,t]\times S)}   \nonumber
 \\
 & \leq &
C_E\left(   \| u_{t}(t_0  ,\cdot  )\|_{  L^2(S) }  + a^{-1}(t_0  )\|   \nabla u (t_0  ,\cdot  ) \|_{  L^2(S) }    +   M (t_0  )   \|    u (t_0  ,\cdot  ) \|_{  L^2(S) }   \right)  \nonumber  \\
&  &+  C_E\| g(x,\tau )  \| _{L^1([t_0  ,t]; L^2(S))}
\end{eqnarray*}
 for all $t > t_0 $, where  $C_E>0$ is a number independent of $t$ and $u$.
\end{proposition}
\medskip

\noindent
{\bf Proof.}
For the energy $E(t) $ (\ref{energy}) of the solution $ u=u(t,s)$ of (\ref{5.7g})
we have
\begin{eqnarray*}
\frac{d}{d t }  E(t )
& = &
 \frac{1}{2}  ( a^{-2}(t) )_t\| \nabla_\sigma  u \|^2_{ L^2( S )}
+  \frac{1}{2}  (M^2 (t) )_t   \|      u   \|^2_{ L^2(S )}
+ \int_{S} (\partial_t u ) g \,d\mu _\sigma  \,.
\end{eqnarray*}
Integration gives
\begin{eqnarray*}
E(t)
& - &
 \int_{t_0}  ^t \Bigg[  \frac{1}{2}  ( a^{-2}(\tau ) )_\tau \| \nabla_\sigma  u \|^2_{ L^2( S )}
+  \frac{1}{2}  (M^2 (\tau ) )_\tau    \|      u   \|^2_{ L^2(S )}   \Bigg] \,d \tau  \\
& = &
E(t_0  ) + \int_{t_0}  ^t\int_{S}    u_t( \tau )  g( \tau )\,d \mu _\sigma \, d t \,.
\end{eqnarray*}
Hence,  we have
\begin{eqnarray*}
 &  &
  \| u_{t}\|^2_{ L^2( S )}+  a^{-2}(t)\| \nabla_\sigma  u \|^2_{ L^2( S )}
 +   M^2(t)   \|      u   \|^2_{ L^2(S )} \\
&  &
- \int_{t_0}  ^t \Bigg[   ( a^{-2}(\tau ) )_\tau \| \nabla_\sigma  u \|^2_{ L^2( S )}
+     (M^2 (\tau ) )_\tau    \|      u   \|^2_{ L^2(S )}    \Bigg] \,d \tau
 \\
& \leq &
  \| u_{t}(t_0  )\|^2_{ L^2( S )}+  a^{-2}(t_0  )\| \nabla_\sigma  u (t_0  )\|^2_{ L^2( S )}
 +  M^2(t_0  )   \|      u (t_0  )  \|^2_{ L^2(S )}\\
 &  &
+2 \max_{t_0  \leq \tau \leq t}  \| u_{t} ( t ) \| _{  L^2(S) }  \int_{t_0} ^t  \| g( \tau )  \| _{  L^2(S) }   \,d\mu _\sigma \, d\tau \,.
\end{eqnarray*}
It follows that for every $\varepsilon >0 $ the following inequality
\begin{eqnarray*}
 &  &
 \| u_{t}\| _{ L^2( S )}+  a^{-1}(t)\| \nabla_\sigma  u \| _{ L^2( S )}
 +   M (t)   \|      u   \| _{ L^2(S )} \\
 &  &
+ \left| \int_{t_0}  ^t \Bigg[    ( a^{-2}(\tau ) )_\tau \| \nabla_\sigma  u \|^2_{ L^2( S )}
+    (M^2 (\tau ) )_\tau    \|      u   \|^2_{ L^2(S )}    \Bigg] \,d \tau\right|^{1/2}\\
& \lesssim  &
  \| u_{t}(t_0  )\| _{ L^2( S )}+  a^{-1}(t_0  )\| \nabla_\sigma  u (t_0  )\| _{ L^2( S )}
 +   M (t_0  )   \|      u (t_0  )  \| _{ L^2(S )}\\
 &  &
+ \varepsilon \max_{t_0  \leq \tau \leq t}  \| u_{t} ( t ) \| _{  L^2(S) }   +  \frac{1}{4\varepsilon } \int_{t_0} ^t  \| g( \tau )  \| _{  L^2(S) }   \,   d\tau 
\end{eqnarray*}
is fulfilled. In fact,  we obtain
\begin{eqnarray*}
 &  &
 \| u_{t}\| _{ L^2( S )}+  a^{-1}(t)\| \nabla_\sigma  u \| _{ L^2( S )}
 +   M (t)   \|      u   \| _{ L^2(S )} \\
 &  &
+ \left(   \|\sqrt{ \dot{a} a^{-3} } \,\nabla_\sigma  u \|^2_{ L^2([t_0  ,t]\times  S )}
+        \|   \sqrt{ |(M^2 (t ) )_t |} \,  u   \|^2_{ L^2([t_0  ,t]\times S )}     \right)^{1/2}\\
& \lesssim &
   \| u_{t}(t_0  )\| _{ L^2( S )}+  a^{-1}(t_0  )\| \nabla_\sigma  u (t_0  )\| _{ L^2( S )}
 +   M (t_0  )   \|      u (t_0  )  \| _{ L^2(S )}\\
 &  &
+   \varepsilon \max_{t_0  \leq \tau \leq t}  \| u_{t} (x,t ) \| _{  L^2(S) }   + \frac{1}{4\varepsilon } \int_{t_0}^t  \| g(x,\tau )  \| _{  L^2(S) }   \,   d\tau  \,,
\end{eqnarray*}
and, consequently, for sufficiently small $\varepsilon >0 $ we have
\begin{eqnarray*}
 &  &
 \| u_{t}\| _{ L^2( S )}+  a^{-1}(t)\| \nabla_\sigma  u \| _{ L^2( S )}
 +   M (t)   \|      u   \| _{ L^2(S )} \\
 &  &
+  \|\sqrt{ \dot{a} a^{-3} }\,\nabla_\sigma  u \|_{ L^2([t_0 ,t]\times  S )}
+        \|   \sqrt{ |\dot{c} | } \,  u   \|_{ L^2([t_0 ,t]\times S )}   \\
& \lesssim  &
  \| u_{t}(t_0 )\| _{ L^2( S )}+  a^{-1}(t_0 )\| \nabla_\sigma  u (t_0 )\| _{ L^2( S )}
 +   M (t_0 )   \|      u (t_0 )  \| _{ L^2(S )}
+    \| g  \| _{  L^1([t_0 ,t];L^2(S)) }    \,.
\end{eqnarray*}
The proposition is proved.  \hfill $\square$

\begin{corollary}
\label{C2.2}
Under condition of the proposition we have
\begin{eqnarray*}
    \| u   \| _{X(t)} +    \| \sqrt{ |\dot{c} | } \, u   \| _{L^2([t_0 ,t]\times S)}
 & \leq &
C_E\Big(   \| u_{t}(t_0 ,\cdot  )\|_{L^2(S)}   + a^{-1}(t_0 )\|   \nabla u (t_0 ,\cdot  ) \|_{L^2(S)}   \\
&  &
 +   M (t_0 )   \|    u (t_0 ,\cdot  ) \|_{L^2(S)}     +   \| g(x,\tau )  \| _{L^1([t_0 ,t]; L^2(S))}  \Big)  \nonumber
\end{eqnarray*}
 for all $t > t_0$.  In particular,
\begin{eqnarray*}
    \| u  \| _{X(t) }
 & \leq &
C_E\Big(   \| u_{t}(t_0 ,\cdot  )\|_{L^2(S)}   + a^{-1}(t_0 )\|   \nabla u (t_0 ,\cdot  ) \|_{L^2(S)}
  +   M (t_0 )   \|    u (t_0 ,\cdot  ) \|_{L^2(S)}    \Big) \nonumber \\
&  &
 +  C_E\| g(x,\tau )  \| _{L^1([t_0 ,t]; L^2(S))}
\end{eqnarray*}
 for all $t > t_0$.
\end{corollary}

\medskip

\section{Estimate of self-interaction term}
\label{S3}

\setcounter{equation}{0}
\renewcommand{\theequation}{\thesection.\arabic{equation}}

 The next lemma is a simple consequence of the  Gagliardo-Nirenberg inequality (see, e.g., \cite{Shatah}, p.22 and Lemma~8.2~\cite{Tartar}) and we give detailed proof for
the sake of the self-completeness of this paper.
\begin{lemma}
\label{L3.1}
\label{L3.1m}
Assume that $\alpha \geq 0 $. If $| F( s, \varphi  ) |\leq C |   \varphi   |^{\alpha +1} $ for all  $(s,\varphi ) \in S\times {\mathbb R} $, then  the inequality
\begin{eqnarray}
\label{G-Nm}
\|  F( s , \varphi  (s ) )\|_{ L^2(S) }
& \leq  &
C_1       \|    \varphi  (s)  \|_{ H^1(S) } ^{\frac{n\alpha }{2}}
  \|    \varphi  (s) \|_{ L^2(S) }  ^{ \alpha +1-\frac{n\alpha }{2 } }
\end{eqnarray}
holds provided that $n \geq 3 $ and  $ 0< \alpha \leq \frac{2}{n-2} $. Moreover, if $F(s ,\varphi ) $ is Lipschitz continuous in $ \varphi $  with exponent $\alpha  $, and $n \geq 3 $, $ 0< \alpha \leq \frac{2}{n-2} $, then
\begin{eqnarray}
\label{G-Ndiffm}
&  &
  \|  F( s , \varphi_1 (s)  )- F( s , \varphi_2 (s)  )\|_{ L^2(S) } \\
& \leq  &
C_1  \max_{\theta  =  \varphi_1, \varphi_2}
\Bigg(   \|    \theta   (s) \|^{\frac{\alpha }{\alpha +1}}_{ H^1(S) }   \,
\|    \varphi_1  (s) -  \varphi_2  (s)  \|^{\frac{1}{\alpha +1}}_{ H^1(S) }
  \Bigg)^{\frac{n\alpha }{2}} \nonumber \\
\hspace{-0.5cm} &  &
\qquad \times\max_{\theta  =  \varphi_1, \varphi_2} \left(  \|  \theta    (s) \|^{\frac{\alpha }{\alpha +1}}_{ L^2(S) } \, \|  \varphi_1  (s) - \varphi_2  (s) \|^{\frac{1}{\alpha +1}}_{ L^2(S) }\right)^{\alpha +1-\frac{n\alpha }{2 }} \,.\nonumber
\end{eqnarray}
\end{lemma}
\medskip 

\noindent
{\bf Proof.} First we consider the case of $S={\mathbb R}^{n} $. Thus, we have to prove that   
 if $| F( x, \varphi  ) |\leq C |   \varphi   |^{\alpha +1} $ for all $(x,\varphi ) \in {\mathbb R}^{n+1}$, then  the inequality
\begin{eqnarray}
\label{G-N}
  \|  F( x, \varphi  (x) )\|_{ L^2({\mathbb R}^n) }
& \leq  &
C_1  \|  \nabla \varphi  (x) \|^{\frac{n\alpha }{2}}_{ L^2({\mathbb R}^n) }  \|  \varphi  (x) \|^{ \alpha +1-\frac{n\alpha }{2 } }_{ L^2({\mathbb R}^n) }
\end{eqnarray}
holds provided that $n \geq 3 $ and  $ 0< \alpha \leq \frac{2}{n-2} $. Moreover, if $F(x,\varphi ) $ is Lipschitz continuous in $ \varphi $ with exponent $\alpha  $, then (\ref{G-Ndiffm}) reads 
\begin{eqnarray}
\label{G-Ndiff}
&  &
\|  F( x, \varphi_1  (x) )- F( x, \varphi_2  (x) )\|_{ L^2({\mathbb R}^n) }\\
& \leq  &
C_1  \max_{\psi =  \varphi_1, \varphi_2}\left( \|  \nabla \psi  (x) \|^{\frac{\alpha }{\alpha +1}}_{ L^2({\mathbb R}^n) } \|  \nabla (\varphi_1  (x) - \varphi_2  (x)) \|^{\frac{1}{\alpha +1}}_{ L^2({\mathbb R}^n) } \right)^{\frac{n\alpha }{2}}  \nonumber \\
&  &
\qquad \times \left(  \|  \psi   (x) \|^{\frac{\alpha }{\alpha +1}}_{ L^2({\mathbb R}^n) } \|  \varphi_1  (x) - \varphi_2  (x) \|^{\frac{1}{\alpha +1}}_{ L^2({\mathbb R}^n) }\right)^{\alpha +1-\frac{n\alpha }{2 }} \,.\nonumber
\end{eqnarray}
\medskip

The proof of the inequality (\ref{G-N}) is very similar to the one of   (\ref{G-Ndiff}) with $\varphi _2=0 $ and it does not require for $F$ to be
Lipschitz continuous in $ \varphi $, therefore  we prove only  the last one.
Let  $1/p_1+1/q_1=1 $, then by  property of $F$ and by H\"older's inequality we have the following inequality
\begin{eqnarray*}
&  &
  \|  F(x, \varphi_1  (x) )-  F(x, \varphi_2  (x) )\|^2_{ L^2({\mathbb R}^n) } \\
  &  \lesssim   &
   \int_{{\mathbb R}^n} |\varphi_1 (x) - \varphi_2  (x)|^2  (|\varphi_1 (x)|+|\varphi_2  (x)|)^{2 \alpha  } \,dx  \\
  &  \lesssim &
  \left(\int_{{\mathbb R}^n} |\varphi_1 (x) - \varphi_2  (x)|^{2p_1}   \,dx \right)^{\frac{1}{p_1}} \left(\int_{{\mathbb R}^n}   ( |\varphi_1 (x)|+|\varphi_2  (x)|)^{2 \alpha q_1 } \,dx \right)^{\frac{1}{q_1}} \,.
 \end{eqnarray*}
Denote $\varphi (x)= \varphi_1 (x) - \varphi_2  (x)$ and $\psi (x)=  |\varphi_1 (x)|+|\varphi_2  (x)|$, then
\[
  \|  F(x, \varphi_1  (x) )-  F(x, \varphi_2  (x) )\| _{ L^2({\mathbb R}^n) }
   \lesssim   
   \|  \varphi   \|_{ L^{2p_1}({\mathbb R}^n)}     \|  \psi \| ^{\alpha } _{ L^{2\alpha q_1}({\mathbb R}^n)} \,.
 \]
 Now we  use the Gagliardo-Nirenberg inequality (see \cite{Shatah}, p.22 and Lemma~8.2~\cite{Tartar})
\[
\|   \varphi  (x) \|_{ L^r({\mathbb R}^n) } \lesssim   \|  \varphi  (x) \|^{1-\vartheta }_{ L^2({\mathbb R}^n) }
 \|  \nabla \varphi  (x) \|^{\vartheta }_{ L^2({\mathbb R}^n) }
\]
 with $0\leq \vartheta \leq 1$.
More precisely, we set
\[
 r_1= 2p_1, \quad \frac{1}{2p_1}  =  \frac{1 }{2}    - \frac{\vartheta_1}{n} ,\qquad r_2= 2\alpha q_1, \quad \frac{1}{2\alpha q_1}  =  \frac{1 }{2}    - \frac{\vartheta_2}{n} \,,
\]
and write
\begin{eqnarray*}
\|   \varphi  (x) \|_{ L^{2p_1} ({\mathbb R}^n) }
& \lesssim  &
\|  \varphi  (x) \|^{1-\vartheta_1 }_{ L^2({\mathbb R}^n) }
 \|  \nabla \varphi  (x) \|^{\vartheta_1 }_{ L^2({\mathbb R}^n) }\,, \\
\|  \psi   (x) \|_{ L^{2\alpha q_1}({\mathbb R}^n) }
& \lesssim  & 
 \|  \psi   (x) \|^{1-\vartheta_2 }_{ L^2({\mathbb R}^n) }
 \|  \nabla \psi   (x) \|^{\vartheta_2 }_{ L^2({\mathbb R}^n) } \,.
\end{eqnarray*}
If we choose
\[
p_1=\alpha +1,\quad q_1= \frac{\alpha +1}{\alpha },\quad \vartheta_1 = \frac{ n \alpha }{ 2(\alpha +1 )}=\vartheta_2\,,
\]
then
\begin{eqnarray*}
\|   \varphi  (x) \|_{ L^{2p_1} ({\mathbb R}^n) }
& \lesssim  &
  \|  \varphi  (x) \|^{1-\frac{ n \alpha }{2(\alpha +1) } }_{ L^2({\mathbb R}^n) }
 \|  \nabla \varphi  (x) \|^{\frac{ n \alpha }{2(\alpha +1)} }_{ L^2({\mathbb R}^n) }  
\end{eqnarray*}
and
\begin{eqnarray*}
\|  \psi   (x) \|_{ L^{2\alpha q_1}({\mathbb R}^n) }
& \lesssim  &
  \|  \psi   (x) \|^{1-\frac{n \alpha  }{2( \alpha +1)}}_{ L^2({\mathbb R}^n) }
 \|  \nabla \psi   (x) \|^{\frac{n \alpha  }{ 2( \alpha +1)}}_{ L^2({\mathbb R}^n) }\,.
\end{eqnarray*}
Hence
\begin{eqnarray*}
&  &
  \|  F(x, \varphi_1  (x) )-  F(x, \varphi_2  (x) )\| _{ L^2({\mathbb R}^n) }  \\
  &  \lesssim  &
  \|  \varphi  (x) \|^{1-\frac{ n \alpha }{2(\alpha +1) } }_{ L^2({\mathbb R}^n) }
 \|  \nabla \varphi  (x) \|^{\frac{ n \alpha }{2(\alpha +1)} }_{ L^2({\mathbb R}^n) }
   \|  \psi   (x) \|^{\alpha(1-  \frac{n \alpha   }{ 2( \alpha +1)})}_{ L^2({\mathbb R}^n) }
 \|  \nabla \psi   (x) \|^{\alpha \frac{n \alpha   }{2( \alpha +1)}}_{ L^2({\mathbb R}^n) }  \,.
 \end{eqnarray*}
Thus, (\ref{G-Ndiff}) is proven.  
\smallskip

 Next we turn to the case of manifold $S$.  Let $U \subseteq {\mathbb R}^n$ be a local chart with local coordinates $x \in U$. If $\psi  \in C_0^\infty (U) $, $\psi  \geq 0 $, then
due to (\ref{G-N}) we have
\begin{eqnarray*}
 \| \psi (x) F( x, \varphi  (x) )\|_{ L^2(U) }
& \lesssim  &
 \|   |\psi (x) \varphi  (x) |^{\alpha +1} \|_{ L^2(U) }\\
& \lesssim  &
  \|   \nabla  ( \psi (x) \varphi  (x) ) \|^{\frac{n\alpha }{2}}_{ L^2(U) } \,\|   \psi (x)  \varphi  (x) \|^{ \alpha +1-\frac{n\alpha }{2 } }_{ L^2(U) } \\
& \lesssim  &
 \|     \varphi  (x)   \|^{\frac{n\alpha }{2}}_{ H^1(U) } \, \|     \varphi  (x) \|^{ \alpha +1-\frac{n\alpha }{2 } }_{ L^2(U) }\,.
\end{eqnarray*}
Let $\{ \psi _j\} $ be a locally finite partition of unity on $S$ subordinated  to the cover $\{U_j\} $. Then due to (\ref{G-N}) we obtain
\begin{eqnarray*}
 \|  F( x, \varphi  (x) )\|_{ L^2(S) }
& =  &
\sum_{j}  \|\psi_j (x) F( x,  \varphi  (x) )\|_{ L^2(U_j) }  \\
& \lesssim  &
 \sum_{j}   \|     \varphi  (x)  \|^{\frac{n\alpha }{2}}_{ H^1(U_j) } \|    \varphi  (x) \|^{ \alpha +1-\frac{n\alpha }{2 } }_{ L^2(U_j) }\,.
\end{eqnarray*}
Since $\frac{2}{n-2 }\geq \alpha  $,   we have $\alpha +1-\frac{n\alpha }{2 }\geq 0$, and, consequently,    due to the inequality $ (\sum_{n} a_n)^k \lesssim \sum_{n} a_n^k \lesssim  (\sum_{n} a_n)^k$
for non-negative numbers $k$, $a_1$,$a_2$,\dots \,,$a_n$ , we obtain
\begin{eqnarray*}
 \|  F( x, \varphi  (x) )\|_{ L^2(S) }
& \lesssim  &
    \left( \sum_{j} \|     \varphi  (s)  \|_{ H^1(U_j) }\right)^{\frac{n\alpha }{2}}
 \left(  \sum_{j}\|    \varphi  (s) \|_{ L^2(U_j) } \right)^{ \alpha +1-\frac{n\alpha }{2 } } \\
& \lesssim  &
      \|     \varphi  (s)  \|_{ H^1(S) } ^{\frac{n\alpha }{2}}
  \|    \varphi  (s) \|_{ L^2(S) }  ^{ \alpha +1-\frac{n\alpha }{2 } } \,.
\end{eqnarray*}
Hence, (\ref{G-Nm}) is proven. Now for $\psi  \in C_0^\infty (U) $, $\psi  \geq 0 $, consider
\begin{eqnarray*}
&  &
  \|   \psi (x)F( x ,\varphi_1  (x)) -  \psi (x)F( x ,\varphi_2 (x)  )  \|_{ L^2(U) } \\
& \lesssim  &
   \max_{\theta  =  \varphi_1, \varphi_2}\left( \|  \nabla (\psi (x) \theta   (x)) \|^{\frac{\alpha }{\alpha +1}}_{ L^2(U) }
\|  \nabla (\psi (x)\varphi_1  (x) - \psi (x)\varphi_2  (x)) \|^{\frac{1}{\alpha +1}}_{ L^2(U) } \right)^{\frac{n\alpha }{2}}  \\
\hspace{-0.5cm} &  &
\qquad \times \left(  \|  \theta    (x) \|^{\frac{\alpha }{\alpha +1}}_{ L^2(U) } \|  \varphi_1  (x) - \varphi_2  (x) \|^{\frac{1}{\alpha +1}}_{ L^2(U) }\right)^{\alpha +1-\frac{n\alpha }{2 }} \,,\nonumber
\end{eqnarray*}
where (\ref{G-Ndiff}) has been used. If $\{\psi  _j\} $ is a locally finite partition of unity on $S$ subordinated  to the cover $\{U_j\} $, then
\begin{eqnarray*}
&  &
 \|   F( s , \varphi_1 (s)  )- F( s , \varphi_2 (s)   )  \|_{ L^2(S) } \\
& = &
 \sum_j \|  \psi_j (x)( F( x ,\varphi_1 (x)  )- F( x , \varphi_2  (x) )) \|_{ L^2(U_j) } \\
& \lesssim  &
   \sum_j
  \left(\max_{\theta  =  \varphi_1, \varphi_2} \|  \nabla  (\psi_j (x)\theta   (x)) \|^{\frac{\alpha }{\alpha +1}}_{ L^2(U_j) }
\|  \nabla (\psi_j (x)\varphi_1  (x) - \psi_j (x)\varphi_2  (x)) \|^{\frac{1}{\alpha +1}}_{ L^2(U_j) } \right)^{\frac{n\alpha }{2}}  \\
\hspace{-0.5cm} &  &
\qquad \times \left( \max_{\theta  =  \varphi_1, \varphi_2} \|  \theta    (x) \|^{\frac{\alpha }{\alpha +1}}_{ L^2(U_j) } \|  \varphi_1  (x) - \varphi_2  (x) \|^{\frac{1}{\alpha +1}}_{ L^2(U_j) }\right)^{\alpha +1-\frac{n\alpha }{2 }}  \\
& \lesssim  &
   \left(\sum_j  \max_{\theta  =  \varphi_1, \varphi_2} \|     \theta   (x) \|^{\frac{\alpha }{\alpha +1}}_{ H^1(U_j) }
\|    (\varphi_1  (x) - \varphi_2  (x)) \|^{\frac{1}{\alpha +1}}_{ H^1(U_j) } \right)^{\frac{n\alpha }{2}}  \\
\hspace{-0.5cm} &  &
\qquad \times \left( \sum_j \max_{\theta  =  \varphi_1, \varphi_2} \|  \theta    (x) \|^{\frac{\alpha }{\alpha +1}}_{ L^2(U_j) } \|  \varphi_1  (x) - \varphi_2  (x) \|^{\frac{1}{\alpha +1}}_{ L^2(U_j) }\right)^{\alpha +1-\frac{n\alpha }{2 }}  \\
& \lesssim  &
   \left(  \max_{\theta  =  \varphi_1, \varphi_2} \|     \theta   (x) \|^{\frac{\alpha }{\alpha +1}}_{ H^1(S) }
\|   \varphi_1  (x) - \varphi_2  (x)  \|^{\frac{1}{\alpha +1}}_{ H^1(S) } \right)^{\frac{n\alpha }{2}}  \\
\hspace{-0.5cm} &  &
\qquad \times \left(  \max_{\theta  =  \varphi_1, \varphi_2} \|  \theta    (x) \|^{\frac{\alpha }{\alpha +1}}_{ L^2(S) } \|  \varphi_1  (x) - \varphi_2  (x) \|^{\frac{1}{\alpha +1}}_{ L^2(S) }\right)^{\alpha +1-\frac{n\alpha }{2 }}
\end{eqnarray*}
proves (\ref{G-Ndiffm}). Lemma is proven. \hfill $\square $
\smallskip

The next proposition gives the estimate of the self-interaction term, which is  transformed   by the reduction of the damping term $n \dot{a} (t)/a (t) \psi _t$ of the equation.
It is also connected with the energy of the linear equation. In the next proposition the norm (\ref{normX}) with $S= {\mathbb R}^n$ will be used.  
Define the   function $\widetilde{C}_{\alpha  ,\Gamma , \alpha_0  } (T)$   as follows
\begin{eqnarray}
\label{CbG}
&  &
\widetilde{C}_{\alpha  ,\Gamma , \alpha_0  } (T)
  =  
\cases{ 1, \quad  \mbox{if} \quad (\ref{GammaB}) \quad  \mbox{and} \quad   \alpha _0=\frac{4}{n}, \cr
 C_{a,\Gamma ,\alpha _0}(T)  
,\quad  \mbox{if} \quad (\ref{GammaB}) \,\,  or \,\, (\ref{GammaU})  \,\, \mbox{and} \quad   \alpha _0<\frac{4}{n}\,.}   
\end{eqnarray}
\begin{proposition}
\label{P3.2}
Assume that $n \alpha \geq 4  $. If $| F( x, \varphi  ) |\leq C |   \varphi   |^{\alpha +1} $ for all $(x,\varphi ) \in {\mathbb R}^{n+1}$, then  the inequality
\begin{eqnarray*} 
 \|a^{\frac{n}{2}} (t) \Gamma (t)F(x, a^{-\frac{n}{2}} (t)   u )  \|_{L^1([t_0 ,T]; L^2({\mathbb R}^n))}   
& \lesssim  &
 \widetilde{C}_{\alpha  ,\Gamma , \alpha_0  } (T)\parallel u \parallel _{X(T)}^{\alpha +1} \,,
\end{eqnarray*}
holds for all $T \in (t_0 ,\infty)$,  with the functions $a(t)$ and $\Gamma (t) $ satisfying conditions (\ref{GammaB})  or   (\ref{GammaU})
of Theorem~\ref{T1.1}.

Moreover, if, additionally,  $F(x,\varphi ) $ is Lipschitz continuous in $ \varphi $ with exponent $\alpha  $, then
\begin{eqnarray}
 &  &
 \|a^{\frac{n}{2}} (t) \Gamma (t)\left( F(x, a^{-\frac{n}{2}} (t)   u ) - F(x, a^{-\frac{n}{2}} (t)   v )\right)   \|_{L^1([t_0 ,T]; L^2({\mathbb R}^n))} \nonumber  \\
& \lesssim  &
\label{3.4}
 \widetilde{C}_{\alpha  ,\Gamma , \alpha_0  } (T)\left( \max_{w=u,v} \parallel w \parallel _{X(t)}^{\alpha } \right) \parallel u -v \parallel _{X(t)}\,.
\end{eqnarray}
\end{proposition}
\medskip

\noindent
{\bf Proof.} We prove only the second part of the proposition since the proof of the first one is very similar, but  it does not appeal to the above mentioned
additional condition on $F$. Denote $b(t)= a^{-\frac{n}{2}} (t) $, $\widetilde{\alpha } := \frac{n\alpha }{2(\alpha +1)}$.
From (\ref{G-Ndiff})  we derive
\begin{eqnarray*}
 &  &
 \|b^{-1} (t) \Gamma (t) \left( F(x, b(t)  u )- F(x, b(t)  v )\right) \|_{L^1([t_0 ,T]; L^2({\mathbb R}^n))}   \\
& \leq  &
\int_{t_0} ^T b^{-1} (t)\Gamma (t) \|  F(x, b(t)  u )- F(x, b(t)  v )  \|_{ L^2({\mathbb R}^n) } \, dt \\
& \lesssim  &
\int_{t_0} ^T b^{\alpha } (t)\Gamma (t) \max_{w =  u,v}\left( \|  \nabla w  (x,t) \|^{\frac{\alpha }{\alpha +1}}_{ L^2({\mathbb R}^n) } \|  \nabla (u  (x,t) - v(x,t)) \|^{\frac{1}{\alpha +1}}_{ L^2({\mathbb R}^n) } \right)^{\widetilde{\alpha }(\alpha +1)}  \nonumber \\
&  &
\qquad \times \left(  \|  w   (x) \|^{\frac{\alpha }{\alpha +1}}_{ L^2({\mathbb R}^n) } \|  u  (x,t) - v(x,t) \|^{\frac{1}{\alpha +1}}_{ L^2({\mathbb R}^n) }\right)^{(1-\widetilde{\alpha })(\alpha +1)}  \, dt \,.
\end{eqnarray*}
For the positive function $\mu =\mu (t) $ this leads to the following estimate
\begin{eqnarray}
\label{3.6}
 &  &
 \|b^{-1} (t) \Gamma (t) \left( F(x, b(t)  u )- F(x, b(t)  v )\right) \|_{L^1([t_0 ,T]; L^2({\mathbb R}^n))}   \\
& \lesssim  &
 \int_{t_0} ^T \mu (t)b^{\alpha } (t)\Gamma (t)  \nonumber \\
&  &
\max_{w =  u,v}\left( \| \mu^{-\frac{2 }{n\alpha  }} (t) \nabla w  (x,t) \|^{\frac{\alpha }{\alpha +1}}_{ L^2({\mathbb R}^n) } \| \mu^{-\frac{2 }{n\alpha  }} (t) \nabla (u  (x,t) - v(x,t)) \|^{\frac{1}{\alpha +1}}_{ L^2({\mathbb R}^n) } \right)^{\widetilde{\alpha }(\alpha +1)}  \nonumber \\
&  &
\qquad \times \left(  \|  w   (x) \|^{\frac{\alpha }{\alpha +1}}_{ L^2({\mathbb R}^n) } \|  u  (x,t) - v(x,t) \|^{\frac{1}{\alpha +1}}_{ L^2({\mathbb R}^n) }\right)^{(1-\widetilde{\alpha })(\alpha +1)}  \, dt\, \nonumber \,,
\end{eqnarray}
where 
\[
\frac{1 }{\widetilde{\alpha }(\alpha +1) } =\frac{2 }{n\alpha  } \,.
\]
Now we set
\begin{eqnarray}
\label{mu}
 \mu (t)  = a(t)  ^{\frac{2n \alpha +n\alpha _0  }{4} }    |\dot{a}(t)| ^{-\frac{ n\alpha _0 }{4} }  
\end{eqnarray}
and consider two cases:  (B)  (\ref{GammaB}) and $  \alpha _0 =  \frac{4}{n}$  ;   (U) $ 0< \alpha _0 <  \frac{4}{n}$ and (\ref{GammaB}) or (\ref{GammaU}).
\smallskip

In the first case (B)  due to (\ref{GammaB}) we obtain
\begin{eqnarray*}
&  &
\mu (t)  b^{\alpha } (t)|\Gamma (t)| \leq const  \,.
\end{eqnarray*}
Therefore,
\begin{eqnarray*}
 &  &
 \|b^{-1} (t) \Gamma (t) \left( F(x, b(t)  u )- F(x, b(t)  v )\right) \|_{L^1([t_0 ,T]; L^2({\mathbb R}^n))}   \\
& \lesssim  &
\int_{t_0} ^T  \max_{w =  u,v}\left( \| \mu^{-\frac{2}{n\alpha   }} (t) \nabla w  (x,t) \|^{\frac{\alpha }{\alpha +1}}_{ L^2({\mathbb R}^n) } \| \mu^{-\frac{2}{n\alpha   }} (t) \nabla (u  (x,t) - v(x,t)) \|^{\frac{1}{\alpha +1}}_{ L^2({\mathbb R}^n) } \right)^{\frac{n\alpha }{2}}  \nonumber \\
&  &
\qquad \times \left(  \|  w   (x) \|^{\frac{\alpha }{\alpha +1}}_{ L^2({\mathbb R}^n) } \|  u  (x,t) - v(x,t) \|^{\frac{1}{\alpha +1}}_{ L^2({\mathbb R}^n) }\right)^{\alpha +1  - \frac{n\alpha }{2 }}  \, dt \\
& \lesssim  &
 \max_{w =  u,v}\left( \|  w   (x) \|^{\frac{\alpha }{\alpha +1}}_{ L^\infty ((t_0 ,T); L^2({\mathbb R}^n)) } \|  u  (x,t) - v(x,t) \|^{\frac{1}{\alpha +1}}_{ L^\infty ((t_0 ,T); L^2({\mathbb R}^n))}\right)^{\alpha +1  - \frac{n\alpha }{2 }} \\
 &  &
 \times \int_{t_0} ^T  \max_{w =  u,v}\left( \| \mu^{-\frac{2}{n\alpha   }} (t) \nabla w  (x,t) \|^{\frac{\alpha }{\alpha +1}}_{ L^2({\mathbb R}^n) } \| \mu^{-\frac{2}{n\alpha   }} (t) \nabla (u  (x,t) - v(x,t)) \|^{\frac{1}{\alpha +1}}_{ L^2({\mathbb R}^n) } \right)^{\frac{n\alpha }{2}}  \, dt  \,.
\end{eqnarray*}
Consider the last integral, where    $\mu (t) $  is replaced with its definition (\ref{mu}). After long but simple transformations  we arrived at 
\begin{eqnarray*}
 &  &
\int_{t_0} ^T  \max_{w =  u,v}\Big( \| \mu^{-\frac{2}{n\alpha   }} (t) \nabla w  (x,t) \|^{\frac{\alpha }{\alpha +1}}_{ L^2({\mathbb R}^n) } \| \mu^{-\frac{2}{n\alpha   }} (t) \nabla (u  (x,t) - v(x,t)) \|^{\frac{1}{\alpha +1}}_{ L^2({\mathbb R}^n) } \Big)^{\frac{n\alpha }{2}}  \, dt   \\
& = &
\int_{t_0} ^T    \max_{w =  u,v}\Big( \|a(t)  ^{-1}\nabla w  (x,t) \|^{\frac{\alpha }{\alpha +1}}_{ L^2({\mathbb R}^n) }
\times  \| a(t)  ^{-1}\nabla (u  (x,t) - v(x,t)) \|^{\frac{1}{\alpha +1}}_{ L^2({\mathbb R}^n) } \Big)^{\frac{n(\alpha-\alpha _0) }{2}}  \\
&  &
\times \max_{w =  u,v}\Big( \|\sqrt{\dot{a}(t) a(t)  ^{-  3 }}\nabla w  (x,t) \|^{\frac{\alpha }{\alpha +1}}_{ L^2({\mathbb R}^n) } \| \sqrt{\dot{a}(t) a(t)  ^{-  3 }}\nabla (u  (x,t) - v(x,t)) \|^{\frac{1}{\alpha +1}}_{ L^2({\mathbb R}^n) } \Big)^{\frac{n\alpha_0 }{2}} \, dt \\
& \lesssim  &
 \max_{w =  u,v}\Big( \|a(t)  ^{-1}\nabla w  (x,t) \|^{\frac{\alpha }{\alpha +1}}_{ L^\infty ((t_0 ,T); L^2({\mathbb R}^n))  }\\
&  &
\hspace{2cm}  \times 
  \| a(t)  ^{-1}\nabla (u  (x,t) - v(x,t)) \|^{\frac{1}{\alpha +1}}_{ L^\infty ((t_0 ,T); L^2({\mathbb R}^n))  } \Big)^{\frac{n(\alpha-\alpha _0) }{2}}  \\
&  &
\qquad \times \int_{t_0} ^T \max_{w =  u,v}\Big( \|\sqrt{\dot{a}(t) a(t)  ^{-  3 }}\nabla w  (x,t) \|^{\frac{\alpha }{\alpha +1}}_{ L^2({\mathbb R}^n) }\\
&  &
\hspace{2cm} \times  \| \sqrt{\dot{a}(t) a(t)  ^{-  3 }}\nabla (u  (x,t) - v(x,t)) \|^{\frac{1}{\alpha +1}}_{ L^2({\mathbb R}^n) } \Big)^{\frac{n\alpha_0 }{2}} \, dt \,.
\end{eqnarray*}
Thus,
\begin{eqnarray*}
 &  &
 \|b^{-1} (t) \Gamma (t) \left( F(x, b(t)  u )- F(x, b(t)  v )\right) \|_{L^1([t_0 ,T]; L^2({\mathbb R}^n))}   \\
& \lesssim  &
 \max_{w =  u,v}\left( \|  w   (x) \|^{\frac{\alpha }{\alpha +1}}_{ L^\infty ((t_0 ,T); L^2({\mathbb R}^n)) } \|  u  (x,t) - v(x,t) \|^{\frac{1}{\alpha +1}}_{ L^\infty ((t_0 ,T); L^2({\mathbb R}^n))}\right)^{\alpha +1  - \frac{n\alpha }{2 }} \\
 &  &
 \times
 \max_{w =  u,v}\Big( \|a(t)  ^{-1}\nabla w  (x,t) \|^{\frac{\alpha }{\alpha +1}}_{ L^\infty ((t_0 ,T); L^2({\mathbb R}^n))  }
  \| a(t)  ^{-1}\nabla (u  (x,t) - v(x,t)) \|^{\frac{1}{\alpha +1}}_{ L^\infty ((t_0 ,T); L^2({\mathbb R}^n))  } \Big)^{\frac{n(\alpha-\alpha _0 )}{2}}  \\
&  &
\times \int_{t_0} ^T \max_{w =  u,v}\Big( \|\sqrt{\dot{a}(t) a(t)  ^{-  3 }}\nabla w  (x,t) \|^{\frac{\alpha }{\alpha +1}}_{ L^2({\mathbb R}^n) } 
\| \sqrt{\dot{a}(t) a(t)  ^{-  3 }}\nabla (u  (x,t) - v(x,t)) \|^{\frac{1}{\alpha +1}}_{ L^2({\mathbb R}^n) } \Big)^{\frac{n\alpha_0 }{2}} \, dt \,.
\end{eqnarray*}
We set here   $  \alpha _0=\frac{4}{n} \leq \alpha $  and by means of  the condition   $\alpha \leq \frac{2}{n-2} $
and by the H\"older inequality with
\[
p_2= \frac{4(\alpha +1)}{n\alpha _0\alpha }\,,\quad q_2= \frac{4(\alpha +1)}{n\alpha _0 }\,,
\]
derive for the integral of the  last inequality the following estimate
\begin{eqnarray*}
 &  &
\int_{t_0} ^T \max_{w =  u,v}\Big( \|\sqrt{\dot{a}(t) a(t)  ^{-  3 }}\nabla w  (x,t) \|^{\frac{\alpha }{\alpha +1}}_{ L^2({\mathbb R}^n) } 
\| \sqrt{\dot{a}(t) a(t)  ^{-  3 }}\nabla (u  (x,t) - v(x,t)) \|^{\frac{1}{\alpha +1}}_{ L^2({\mathbb R}^n) } \Big)^{\frac{n\alpha_0 }{2}} \, dt  \\
&  & \leq 
\max_{w =  u,v} \left(     \|\sqrt{\dot{a}(t) a(t)  ^{-  3 }}\nabla w  (x,t) \| _{ L^2([t_0 ,T]\times {\mathbb R}^n) }  \, dt \right)^{\frac{n \alpha _0\alpha }{2(\alpha +1) }} \\
&  &
\quad \times \left(  \| \sqrt{\dot{a}(t) a(t)  ^{-  3 }}\nabla (u  (x,t) - v(x,t)) \| _{ L^2([t_0 ,T]\times {\mathbb R}^n) }  \, dt  \right) ^{\frac{n \alpha _0  }{2(\alpha +1) }}\,.
\end{eqnarray*}
This implies 
\begin{eqnarray*}
\hspace{-0.5cm}  &  &
 \|a^{\frac{n}{2}} (t) \Gamma (t)\left( F(x, a^{-\frac{n}{2}} (t)   u ) - F(x, a^{-\frac{n}{2}} (t)   v )\right)   \|_{L^1([t_0 ,T]; L^2({\mathbb R}^n))} \nonumber  \\
\hspace{-0.5cm} & \lesssim  &
\label{3.3}
\widetilde{C}_{\alpha  ,\Gamma , \alpha_0  } (T)
\max_{w=u,v}  \left(   \|  w (x,t) \|_{ L^\infty((t_0 ,T);  L^2({\mathbb R}^n)) }  ^{ \frac{\alpha }{\alpha +1}}  \|   u  -v  \|^{\frac{1}{\alpha +1}}_{ L^\infty((t_0 ,T);  L^2({\mathbb R}^n)) } \right)^{ \alpha +1-\frac{n\alpha }{2} }  \nonumber    \\
\hspace{-0.5cm} &  &
\times   \left(   \| a^{-1} (t)  \nabla_x  w (x,t) \|^{ \frac{\alpha }{\alpha +1} } _{ L^\infty((t_0 ,T);  L^2({\mathbb R}^n)) }
 \| a^{-1} (t)  \nabla_x  (u (x,t) -v(x,t))\|^{  \frac{1}{\alpha +1} }_{ L^\infty((t_0 ,T);  L^2({\mathbb R}^n)) }  \right)^{ \frac{n(\alpha -\alpha _0 )}{2} }  \nonumber   \\
\hspace{-0.5cm} &  &
\times   \left(
  \|    \sqrt{ \dot{a}  a^{-3}  }  \nabla_x  w (x,t) \|_{L^2((t_0 ,T)\times {\mathbb R}^n) }  ^{ \frac{\alpha }{\alpha +1} }
\|    \sqrt{ \dot{a}  a^{-3}  }  \nabla_x  (u (x,t) -v(x,t))\|_{L^2((t_0 ,T)\times {\mathbb R}^n) }  ^{  \frac{1}{\alpha +1} }    \right) ^{\frac{ n\alpha _0 }{2} } 
\end{eqnarray*}
and, consequently, (\ref{3.4}).
\smallskip

In the second case (U) with $\alpha_0 < \frac{4}{n}  $   we use  (\ref{3.6})
with $\frac{1}{p_1} + \frac{1}{q_1}=1$,     and   $\widetilde{\alpha }=\frac{n\alpha }{2(\alpha +1)} $.     We choose $\alpha _0, p_1, q_1 $ such that
\[
\frac{1}{p_1} =  1-\frac{n\alpha _0}{4}>0\,,\quad \frac{1}{q_1} =   \frac{n\alpha _0}{4}>0\,,\quad
0<\alpha _0< \frac{4}{n}\,,
\]
and then apply H\"older inequality:
\begin{eqnarray*}
 &  &
 \|b^{-1} (t) \Gamma (t) \left( F(x, b(t)  u )- F(x, b(t)  v )\right) \|_{L^1([t_0 ,T]; L^2({\mathbb R}^n))}   \\
& \lesssim  &
\Bigg( \int_{t_0} ^T [\mu (t)  b^{\alpha } (t) \Gamma (t)]^{p_1} \,dt \Bigg)^{1/p_1} \\
&  &
\times \Bigg( \int_{t_0} ^T \max_{w =  u,v}\left( \| \mu^{-\frac{2 }{n\alpha  }} (t) \nabla w  (x,t) \|^{\frac{\alpha }{\alpha +1}}_{ L^2({\mathbb R}^n) } \| \mu^{-\frac{2 }{n\alpha  }} (t) \nabla (u  (x,t) - v(x,t)) \|^{\frac{1}{\alpha +1}}_{ L^2({\mathbb R}^n) } \right)^{\widetilde{\alpha }(\alpha +1)q_1}   \\
&  &
\qquad \times \left(  \|  w   (x) \|^{\frac{\alpha }{\alpha +1}}_{ L^2({\mathbb R}^n) } \|  u  (x,t) - v(x,t) \|^{\frac{1}{\alpha +1}}_{ L^2({\mathbb R}^n) }\right)^{(1-\widetilde{\alpha })(\alpha +1)q_1}  \, dt \Bigg)^{1/q_1}  \,.
\end{eqnarray*}
Next we use definition of the function $\widetilde{C}_{\alpha  ,\Gamma , \alpha_0  } (T) $  (\ref{CbG}) that, in this case,  reads
\begin{eqnarray*}
&    &
\widetilde{C}_{\alpha  ,\Gamma , \alpha_0  } (T) = \left( \int_{t_0} ^T \left[\mu (t)  b^{\alpha } (t)\Gamma (t)\right]^{p_1} \,dt \right)^{\frac{1}{p_1}}  \,.
\end{eqnarray*}
Hence we obtain the estimate
\begin{eqnarray*}
 &  &
 \|b^{-1} (t) \Gamma (t) \left( F(x, b(t)  u )- F(x, b(t)  v )\right) \|_{L^1([t_0 ,T]; L^2({\mathbb R}^n))}   \\
& \lesssim  &
\widetilde{C}_{\alpha  ,\Gamma , \alpha_0  } (T) \\
&  &
\times \Bigg( \int_{t_0} ^T \max_{w =  u,v}\left( \| \mu^{-\frac{2 }{n\alpha  }} (t) \nabla w  (x,t) \|^{\frac{\alpha }{\alpha +1}}_{ L^2({\mathbb R}^n) } 
\| \mu^{-\frac{2 }{n\alpha  }} (t) \nabla (u  (x,t) - v(x,t)) \|^{\frac{1}{\alpha +1}}_{ L^2({\mathbb R}^n) } \right)^{\widetilde{\alpha }(\alpha +1)q_1}   \\
&  &
\qquad \times \left(  \|  w   (x) \|^{\frac{\alpha }{\alpha +1}}_{ L^2({\mathbb R}^n) } \|  u  (x,t) - v(x,t) \|^{\frac{1}{\alpha +1}}_{ L^2({\mathbb R}^n) }\right)^{(1-\widetilde{\alpha })(\alpha +1)q_1}  \, dt \Bigg)^{1/q_1}  \,.
\end{eqnarray*}
It follows
\begin{eqnarray*}
 &  &
 \|b^{-1} (t) \Gamma (t) \left( F(x, b(t)  u )- F(x, b(t)  v )\right) \|_{L^1([t_0 ,T]; L^2({\mathbb R}^n))}   \\
& \lesssim  &
\widetilde{C}_{\alpha  ,\Gamma , \alpha_0  } (T) \\
&  &
\times \max_{w =  u,v}\left(  \|  w   (x) \|^{\frac{\alpha }{\alpha +1}}_{ L^\infty((t_0 ,T);  L^2({\mathbb R}^n)) } \|  u  (x,t) - v(x,t) \|^{\frac{1}{\alpha +1}}_{ L^\infty((t_0 ,T);  L^2({\mathbb R}^n)) }\right)^{(1-\widetilde{\alpha })(\alpha +1)} \\
&  &
\times \Bigg( \int_{t_0} ^T \max_{w =  u,v}\Big( \| \mu^{-\frac{2 }{n\alpha  }} (t) \nabla w  (x,t) \|^{\frac{\alpha }{\alpha +1}}_{ L^2({\mathbb R}^n) } \\
&  &
\hspace{2cm}\times \| \mu^{-\frac{2 }{n\alpha  }} (t) \nabla (u  (x,t) - v(x,t)) \|^{\frac{1}{\alpha +1}}_{ L^2({\mathbb R}^n) } \Big)^{\widetilde{\alpha }(\alpha +1)q_1}    \, dt \Bigg)^{1/q_1}  \,.
\end{eqnarray*}
Next we estimate the integral of the above inequality:
\begin{eqnarray*}
 &  &
\int_{t_0} ^T \max_{w =  u,v}\left( \| \mu^{-\frac{2 }{n\alpha  }} (t) \nabla w  (x,t) \|^{\frac{\alpha }{\alpha +1}}_{ L^2({\mathbb R}^n) } 
\| \mu^{-\frac{2 }{n\alpha  }} (t) \nabla (u  (x,t) - v(x,t)) \|^{\frac{1}{\alpha +1}}_{ L^2({\mathbb R}^n) } \right)^{\widetilde{\alpha }(\alpha +1)q_1}    \, dt  \\
&  \leq &
\max_{w =  u,v} \Bigg( \| a^{-1}(t) \nabla w  (x,t) \|^{\frac{\alpha }{\alpha +1}}_{ L^\infty((t_0 ,T);  L^2({\mathbb R}^n))} \\
&  &
\hspace{2cm} \times \| a^{-1}(t) \nabla (u  (x,t) - v(x,t)) \|^{\frac{1}{\alpha +1}}_{L^\infty((t_0 ,T);  L^2({\mathbb R}^n)) } \Bigg)^{\frac{n(\alpha-\alpha _0) }{2}q_1}  \\
&  &
\hspace{3cm} \times \int_{t_0} ^T \max_{w =  u,v}\Bigg( \| \sqrt{\dot{a}(t) a(t)  ^{-  3 }} \nabla w  (x,t) \|^{\frac{\alpha }{\alpha +1}}_{ L^2({\mathbb R}^n) } \\
&  &
\hspace{4cm} \times \| \sqrt{\dot{a}(t) a(t)  ^{-  3 }} \nabla (u  (x,t) - v(x,t)) \|^{\frac{1}{\alpha +1}}_{ L^2({\mathbb R}^n) } \Bigg)^{\frac{n\alpha_0}{2}q_1}    \, dt \,.
\end{eqnarray*}
It remains to estimate the last integral, and we set $\alpha _0 =\frac{4}{nq_1} <\alpha  $ and apply H\"older's inequality with
\[
p_2= \frac{4(\alpha +1)}{\alpha _0n\alpha q_1}\,,\quad q_2= \frac{4(\alpha +1)}{\alpha _0n q_1} \,,\quad \alpha _0 =\frac{4}{nq_1} <\alpha \,.
\]
Thus we obtain
\begin{eqnarray*}
 &  &
\int_{t_0} ^T \max_{w =  u,v}\Bigg( \| \sqrt{\dot{a}(t) a(t)  ^{-  3 }} \nabla w  (x,t) \|^{\frac{\alpha }{\alpha +1}}_{ L^2({\mathbb R}^n) } \\
&  &
\hspace{2cm} \times \| \sqrt{\dot{a}(t) a(t)  ^{-  3 }} \nabla (u  (x,t) - v(x,t)) \|^{\frac{1}{\alpha +1}}_{ L^2({\mathbb R}^n) } \Bigg)^{\frac{n\alpha_0}{2}q_1}    \, dt  \\
&  \leq &
\max_{w =  u,v} \left(     \|\sqrt{\dot{a}(t) a(t)  ^{-  3 }}\nabla w  (x,t) \| _{ L^2([t_0 ,T]\times {\mathbb R}^n) }   \right)^{\frac{n \alpha _0\alpha }{2(\alpha +1) }q_1} \\
&  &
\hspace{2cm} \times \left(  \| \sqrt{\dot{a}(t) a(t)  ^{-  3 }}\nabla (u  (x,t) - v(x,t)) \| _{ L^2([t_0 ,T]\times {\mathbb R}^n) }     \right) ^{\frac{n \alpha _0  }{2(\alpha +1) }q_1}\,.
\end{eqnarray*}
To check   (\ref{3.4}) of the statement we just apply the inequality from   Problem 78\cite{Polya-Szego}.
Proposition is proven. \hfill $\square$

\medskip
 
In the special case of the de Sitter spacetime, $a(t)=\exp(Ht)$  and $\Gamma =const$,    the last proposition  implies results of Lemma~3.1~\cite{Nakamura}. The next proposition is the analogue of Proposition~\ref{P3.2}  
but for general manifolds.

\begin{proposition}
\label{P3.5}
Assume that $n \alpha \geq 4    $.
If $| F( s, \varphi  ) |\leq C |   \varphi   |^{\alpha +1} $ for all  $(s,\varphi ) \in S \times  {\mathbb R} $, then  the inequality
\begin{eqnarray}  
 \|a^{\frac{n}{2}} (t) \Gamma (t)F(s , a^{-\frac{n}{2}} (t)   u )  \|_{L^1([t_0 ,T]; L^2(S))}  
& \leq  &
\label{3.14a}
C_N\widetilde{C}_{\alpha  ,\Gamma , \alpha_0  } (T)\parallel u \parallel _{X(T)}^{\alpha +1} \,,
\end{eqnarray}
holds for all $T \in (t_0 ,\infty)$,  with the functions $a(t)$ and $\Gamma (t) $ satisfying conditions (\ref{GammaB})  or   (\ref{GammaU})
of Theorem~\ref{T1.1}. 
   The   function $\widetilde{C}_{\alpha  ,\Gamma , \alpha_0  } (T)$  is defined in (\ref{CbG})

Moreover, if, additionally,  $F(s ,\varphi ) $ is Lipschitz continuous in $ \varphi $ with exponent $\alpha  $, then
\begin{eqnarray}
 &  &
 \|a^{\frac{n}{2}} (t) \Gamma (t)\left( F(s , a^{-\frac{n}{2}} (t)   u ) - F(s , a^{-\frac{n}{2}} (t)   v )\right)   \|_{L^1([t_0 ,T]; L^2(S))} \nonumber  \\
& \leq  &
\label{3.14b}
C_N   \widetilde{C}_{\alpha  ,\Gamma , \alpha_0  } (T)\left(\max_{w=u,v} \parallel w \parallel _{X(t)}^{\alpha } \right) \parallel u -v \parallel _{X(t)}\,.
\end{eqnarray}
\end{proposition}
\medskip

\noindent
{\bf Proof.} Let $U \subseteq {\mathbb R}^n$ be a local chart with local coordinates $x \in U$. If $\varphi \in C_0^\infty (U) $, then
due to  Proposition~\ref{P3.2} we have
\begin{eqnarray*}
 &  &
 \| a^{\frac{n}{2}} (t) \Gamma (t) \varphi (x) F(x, a^{-\frac{n}{2}} (t) u (x.t))  \|_{L^1([t_0 ,T]; L^2({\mathbb R}^n))}  \nonumber  \\
& \lesssim  &
\widetilde{C}_{\alpha  ,\Gamma , \alpha_0  } (T)
 \|  \varphi (x)  u (x,t) \|^{ \alpha +1-\frac{n\alpha }{2} }_{ L^\infty((t_0 ,T);  L^2(U)) }   \nonumber  \\
&  &
\times     \| a^{-1} (t)  \nabla_x (\varphi (x)  u (x,t)) \|^{ \frac{n(\alpha -\alpha _0 )}{2} }_{ L^\infty((t_0 ,T);  L^2(U)) }
  \|    \sqrt{ \dot{a}  a^{-3}  }  \nabla_x (\varphi (x)  u (x,t)) \|_{L^2((t_0 ,T)\times U) }  ^{\frac{ n\alpha _0 }{2} } \\
& \lesssim  &
\widetilde{C}_{\alpha  ,\Gamma , \alpha_0  } (T)
 \|   u (x,t) \|^{ \alpha +1-\frac{n\alpha }{2} }_{ L^\infty((t_0 ,T);  L^2(U)) }       \| a^{-1} (t)    u (x,t)  \|^{ \frac{n(\alpha -\alpha _0 )}{2} }_{ L^\infty((t_0 ,T);  H^1(U)) }\\
&  &
\times \left(   \|    \sqrt{ \dot{a}  a^{-3}  }    u (x,t)  \|_{L^2((t_0 ,T)\times U) }  ^{\frac{ n\alpha _0 }{2} } +\|    \sqrt{ \dot{a}  a^{-3}  }  \nabla_x    u (x,t)  \|_{L^2((t_0 ,T)\times U) }  ^{\frac{ n\alpha _0 }{2} } \right)\,.
\end{eqnarray*}
Let $\{ \psi _j\} $ be a locally finite partition of unity on $S$ subordinated  to the cover $\{U_j\} $. Then
\begin{eqnarray*}
 &  &
 \|  a^{\frac{n}{2}} (t) \Gamma (t)F(x, a^{-\frac{n}{2}} (t)   u )  \|_{L^1([t_0 ,T]; L^2(S))}  \\
& \lesssim  &
\sum_j \|a^{\frac{n}{2}} (t) \Gamma (t) \psi _j F(x, a^{-\frac{n}{2}} (t)  u )  \|_{L^1([t_0 ,T]; L^2(U_j))}  \\
& \lesssim  &
 \widetilde{C}_{\alpha  ,\Gamma , \alpha_0  } (T)
\sum_j  \|   u (x,t) \|^{ \alpha +1-\frac{n\alpha }{2} }_{ L^\infty((t_0 ,T);  L^2(U_j)) }    \| a^{-1} (t)     u (x,t) \|^{ \frac{n(\alpha -\alpha _0 )}{2} }_{ L^\infty((t_0 ,T);  H^1(U_j)) } \\
&  &
\times   
 \left(  \|    \sqrt{ |\dot{a}|  a^{-3}  }     u (x,t) \|_{L^2((t_0 ,T)\times U_j) }  ^{\frac{ n\alpha _0 }{2} }  + \|    \sqrt{ |\dot{a}|  a^{-3}  }  \nabla_x  u (x,t) \|_{L^2((t_0 ,T)\times U_j) }  ^{\frac{ n\alpha _0 }{2} } \right)\\
& \lesssim  &
 \widetilde{C}_{\alpha  ,\Gamma , \alpha_0  } (T)
\left(  \sum_j \|   u (x,t) \|_{ L^\infty((t_0 ,T);  L^2(U_j)) } \right)^{ \alpha +1-\frac{n\alpha }{2} }  \\
&  &
\times  \left(    \sum_j  \| a^{-1} (t)     u (x,t) \|_{ L^\infty((t_0 ,T);  H^1(U_j)) } \right)^{ \frac{n(\alpha -\alpha _0 )}{2} }\\
&  &
\times \left(   \sum_j  \left( \|    \sqrt{ |\dot{a}|  a^{-3}  }     u (x,t) \|_{L^2((t_0 ,T)\times U_j) } +\|    \sqrt{ |\dot{a}|  a^{-3}  }  \nabla_x  u (x,t) \|_{L^2((t_0 ,T)\times U_j) } \right) \right)^{\frac{ n\alpha _0 }{2} }
\,.
\end{eqnarray*}
This proves the inequality
\begin{eqnarray*} 
 &  &
 \|a^{\frac{n}{2}} (t) \Gamma (t)F(s , a^{-\frac{n}{2}} (t)   u )  \|_{L^1([t_0 ,T]; L^2(S))}  \nonumber  \\
& \lesssim  &
 \widetilde{C}_{\alpha  ,\Gamma , \alpha_0  } (T)
 \|   u   \|^{ \alpha +1-\frac{n\alpha }{2} }_{ L^\infty((t_0 ,T);  L^2(S)) }   \nonumber  \\
&  &
\times     \| a^{-1} (t)     u   \|^{ \frac{n(\alpha -\alpha _0 )}{2} }_{ L^\infty((t_0 ,T);  H^1(S)) }
 \left(  \|    \sqrt{ \dot{a}  a^{-3}  }    u   \|_{L^2((t_0 ,T)\times S) }  ^{\frac{ n\alpha _0 }{2} }+  \|    \sqrt{ \dot{a}  a^{-3}  }  \nabla_x  u   \|_{L^2((t_0 ,T)\times S) }  ^{\frac{ n\alpha _0 }{2} }\right) \,.
\end{eqnarray*} 
In order to prove  (\ref{3.14a}) we note that
\begin{eqnarray*}
 \|    \sqrt{  \dot{a}   a^{-3}  }     u (x,t) \| _{L^2((t_0 ,T)\times U_j)} 
& \lesssim  &
 a^{-1}(t_0) (\| M(\cdot )u   \|_{L^\infty((t_0,T);L^2( U_j))} )   \,. 
\end{eqnarray*}
To prove  (\ref{3.14b}) we apply Proposition~\ref{P3.2}  
\begin{eqnarray*}
 &  &
 \|a^{\frac{n}{2}} (t) \Gamma (t)\left( F(s , a^{-\frac{n}{2}} (t)   u ) - F(s , a^{-\frac{n}{2}} (t)   v )\right)   \|_{L^1([t_0 ,T]; L^2(S))} \nonumber  \\
& \lesssim  &
\sum_j  \| a^{\frac{n}{2}} (t) \Gamma (t) \psi _j \left( F(s , a^{-\frac{n}{2}} (t) u ) - F(s , a^{-\frac{n}{2}} (t)   v )\right)   \|_{L^1([t_0 ,T]; L^2(U_j))} \nonumber  \\
& \lesssim  &
 \widetilde{C}_{\alpha  ,\Gamma , \alpha_0  } (T)
\sum_j  \max_{w=u,v}  \left(   \|  w   \|_{ L^\infty((t_0 ,T);  L^2(U_j)) }  ^{ \frac{\alpha }{\alpha +1}}  \|   u  -v  \|^{\frac{1}{\alpha +1}}_{ L^\infty((t_0 ,T);  L^2(U_j)) } \right)^{ \alpha +1-\frac{n\alpha }{2} }  \nonumber    \\
&  &
\times   \left(   \| a^{-1} (t)  \nabla_\sigma  (\psi _j  w ) \|^{ \frac{\alpha }{\alpha +1} } _{ L^\infty((t_0 ,T);  L^2(U_j)) }
 \| a^{-1} (t)  \nabla_\sigma  (\psi _j  (u -v)) \|^{  \frac{1}{\alpha +1} }_{ L^\infty((t_0 ,T);  L^2(U_j)) }  \right)^{ \frac{n(\alpha -\alpha _0 )}{2} }  \nonumber   \\
&  &
\times  \left(
  \|    \sqrt{ \dot{a}  a^{-3}  }  \nabla_\sigma (\psi _j   w ) \|_{L^2((t_0 ,T)\times U_j) }  ^{ \frac{\alpha }{\alpha +1} }
\|    \sqrt{ \dot{a}  a^{-3}  }  \nabla_\sigma  (\psi _j  (u -v)) \|_{L^2((t_0 ,T)\times U_j) }  ^{  \frac{1}{\alpha +1} }    \right) ^{\frac{ n\alpha _0 }{2} }\\
& \lesssim  &
 \widetilde{C}_{\alpha  ,\Gamma , \alpha_0  } (T)
\sum_j  \max_{w=u,v}  \left(   \|  w   \|_{ L^\infty((t_0 ,T);  L^2(U_j)) }  ^{ \frac{\alpha }{\alpha +1}}  \|   u  -v  \|^{\frac{1}{\alpha +1}}_{ L^\infty((t_0 ,T);  L^2(U_j)) } \right)^{ \alpha +1-\frac{n\alpha }{2} }  \nonumber    \\
&  &
\times   \max_{w=u,v}  \left(   \| a^{-1} (t)      w  \|^{ \frac{\alpha }{\alpha +1} } _{ L^\infty((t_0 ,T);  H^1(U_j)) }
 \| a^{-1} (t)      (u -v)  \|^{  \frac{1}{\alpha +1} }_{ L^\infty((t_0 ,T);  H^1(U_j)) }  \right)^{ \frac{n(\alpha -\alpha _0 )}{2} }  \nonumber   \\
&  &
\times  \max_{w=u,v}  \Bigg(
  \|    \sqrt{ \dot{a}  a^{-3}  }      w  \|_{L^2((t_0 ,T)\times U_j) }  ^{ \frac{\alpha }{\alpha +1} }
\|    \sqrt{ \dot{a}  a^{-3}  }     (u -v) \|_{L^2((t_0 ,T)\times U_j) }  ^{  \frac{1}{\alpha +1} }  \\
&  &
+
  \|    \sqrt{ \dot{a}  a^{-3}  }     w  \|_{L^2((t_0 ,T)\times U_j) }  ^{ \frac{\alpha }{\alpha +1} }
\|    \sqrt{ \dot{a}  a^{-3}  }  \nabla_\sigma  (u -v) \|_{L^2((t_0 ,T)\times U_j) }  ^{  \frac{1}{\alpha +1} }  \\
&  &
+  \|    \sqrt{ \dot{a}  a^{-3}  }  \nabla_\sigma   w  \|_{L^2((t_0 ,T)\times U_j) }  ^{ \frac{\alpha }{\alpha +1} }
\|    \sqrt{ \dot{a}  a^{-3}  }    (u -v) \|_{L^2((t_0 ,T)\times U_j) }  ^{  \frac{1}{\alpha +1} }  \\
&  &
+   \|    \sqrt{ \dot{a}  a^{-3}  }  \nabla_\sigma   w  \|_{L^2((t_0 ,T)\times U_j) }  ^{ \frac{\alpha }{\alpha +1} }
\|    \sqrt{ \dot{a}  a^{-3}  }  \nabla_\sigma  (u -v) \|_{L^2((t_0 ,T)\times U_j) }  ^{  \frac{1}{\alpha +1} }   \Bigg) ^{\frac{ n\alpha _0 }{2} }\,.
\end{eqnarray*}
Thus,
\begin{eqnarray*}
 &  &
 \|a^{\frac{n}{2}} (t) \Gamma (t)\left( F(s , a^{-\frac{n}{2}} (t)   u ) - F(s , a^{-\frac{n}{2}} (t)   v )\right)   \|_{L^1([t_0 ,T]; L^2(S))} \nonumber  \\
& \lesssim  &
\label{3.13b}
\widetilde{C}_{\alpha  ,\Gamma , \alpha_0  } (T)
\left(  \left(  \max_{w=u,v}  \|  w   \|_{ L^\infty((t_0 ,T);  L^2(S)) }  ^{ \frac{\alpha }{\alpha +1}}  \right) \|   u  -v  \|^{\frac{1}{\alpha +1}}_{ L^\infty((t_0 ,T);  L^2(S)) } \right)^{ \alpha +1-\frac{n\alpha }{2} } \nonumber    \\
&  &
\times \left( \left(   \max_{w=u,v}   \| a^{-1} (t)      w  \|^{ \frac{\alpha }{\alpha +1} } _{ L^\infty((t_0 ,T);  H^1(S)) }\right)
 \| a^{-1} (t)     (u -v) \|^{  \frac{1}{\alpha +1} }_{ L^\infty((t_0 ,T);  H^1(S)) }  \right)^{ \frac{n(\alpha -\alpha _0 )}{2} }   \nonumber   \\
&  &
\times   \Bigg[   \max_{w=u,v} \Bigg(
  \|    \sqrt{ \dot{a}  a^{-3}  }      w  \|_{L^2((t_0 ,T)\times S) }  ^{ \frac{\alpha }{\alpha +1} } 
\|    \sqrt{\dot{a}  a^{-3}  }     (u -v)\|_{L^2((t_0 ,T)\times S) }  ^{  \frac{1}{\alpha +1} }   \nonumber    \\
&  &
+
  \|    \sqrt{ \dot{a}  a^{-3}  }      w  \|_{L^2((t_0 ,T)\times S) }  ^{ \frac{\alpha }{\alpha +1} }
\|    \sqrt{\dot{a}  a^{-3}  }  \nabla_\sigma   (u -v)\|_{L^2((t_0 ,T)\times S) }  ^{  \frac{1}{\alpha +1} }  \nonumber   \\
&  &
+
  \|    \sqrt{ \dot{a}  a^{-3}  }  \nabla_\sigma   w  \|_{L^2((t_0 ,T)\times S) }  ^{ \frac{\alpha }{\alpha +1} }
\|    \sqrt{\dot{a}  a^{-3}  }      (u -v)\|_{L^2((t_0 ,T)\times S) }  ^{  \frac{1}{\alpha +1} }  \nonumber  \\
&  &
+
  \|    \sqrt{ \dot{a}  a^{-3}  }  \nabla_\sigma   w  \|_{L^2((t_0 ,T)\times S) }  ^{ \frac{\alpha }{\alpha +1} }
\|    \sqrt{\dot{a}  a^{-3}  }  \nabla_\sigma   (u -v)\|_{L^2((t_0 ,T)\times S) }  ^{  \frac{1}{\alpha +1} }   \Bigg)\Bigg] ^{\frac{ n\alpha _0 }{2} }\,.
\end{eqnarray*}
This completes the proof of  (\ref{3.14b}). Proposition is proven.
\hfill $\square$

\bigskip

\section{Completion of the proof of theorems. Examples}
\label{S4}
\setcounter{equation}{0}

\subsection{Integral equation. Proof of theorems.}
 
Now we consider the Cauchy problem for the equation
\begin{eqnarray*}
 &  &
  \square_g \psi
 = m^2  \psi   + V'_\psi (t,x,\psi ) \,.
\end{eqnarray*}
It can be   written as an equation for $u= a^{\frac{n}{2}}(t)  \psi   $:
\begin{eqnarray}
\label{eq4.1}
 &  &
u_{tt}-  \frac{1}{a^2(t)}\Delta _\sigma  u
+M^2(t) u  =  - \frac{1}{b(t)} V'_\psi (t,x, b(t) u  ) \,,
\end{eqnarray}
with $b(t)=a^{-\frac{n}{2}}(t) $ and with the curved mass (\ref{M}).
Here $ \Delta _\sigma $ is Laplace-Beltrami operator in the metric $\sigma $.

Denote by $G$ the solution operator for the Cauchy problem 
\begin{eqnarray*} 
 &  &
u_{tt}-  \frac{1}{a^2(t)}\Delta _\sigma  u
+M^2(t) u  =  f \,,\quad u(t_0,s)=u_t(t_0,s)=0\,.
\end{eqnarray*}
By applying the  operator $G$ to  the  equation  (\ref{eq4.1}), the problem can be rewritten as an  integral equation  for the function $\Phi =\Phi (t,s) $:
\begin{eqnarray}
\label{IE}
  \Phi (t,s)
&  = &
\Phi _0(t,s) -
G \big[  b^{-1} (t) V'_\psi (t,x, b(t) \Phi   )  \big] (t,s)    \,.
\end{eqnarray}
Here $\Phi _0(t,s) $ is a given function.
For the numbers $R>0 $ and $T>t_0  $ we define the complete metric space
\begin{eqnarray*}
 &  &
X({R,T})  := \{ \Phi \in C ([t_0 ,T] ; H_{(1)} (S) )\cap C^1([t_0 ,T] ; L^{2} (S)) \; | \;
 \parallel  \Phi  \parallel _{X(T)}
\le R \}
\end{eqnarray*}
with the metric
\[
d(\Phi _1,\Phi _2) :=  \parallel  \Phi _1   - \Phi _2   \parallel _{X(T)}    \,.
\]

Denote by $ C^{(-1)}_{\mu ,a,\Gamma }(r) $ the function inverse to $ C_{\mu ,a,\Gamma }(T ) $, which is given by the second case of (\ref{CbG}), then $ C^{(-1)}_{\mu ,a,\Gamma }(0) =t_0  $.
The following theorem guarantees  local and global solvability of the   integral equation (\ref{IE}).

\begin{theorem}
\label{T4.1}
(i) Assume that conditions of Theorem~\ref{T1.2} are satisfied. Then
for every $\Phi _0 $ $ \in X(R ,T)$ there exist  $T_1>t_0 $, $R_1>0 $, and the unique (local) solution $\Phi \in X(R_1 ,T)$ of the equation
(\ref{IE}). The life span $T_1 -t_0 >0$ can be estimated from below as follows: there is $ C  $ such that for every $R_1 >R$,
\[
T_1 -t_0  \geq C  \min \left\{ C^{(-1)}_{\mu ,a,\Gamma }\left(\frac{R_1-R}{c_0 R} \right) ,
C^{(-1)}_{\mu ,a,\Gamma } \left( \frac{ 1}{c_0 R_1^\alpha }  \right)   \right\}.
\]

(ii) Assume that conditions of Theorem~\ref{T1.1} are satisfied, then there is  $\varepsilon _0>0 $ such that for every given function $\Phi _0 \in X(\varepsilon ,T)$ with small norm
\[
\parallel \Phi _0 \parallel _{X(T)} \leq \varepsilon  < \varepsilon _0 ,
\]
$ 0< T \leq \infty$, the integral equation (\ref{IE}) has a unique solution $\Phi  \in X(2\varepsilon ,T)$ and
\[
\parallel \Phi   \parallel _{X(T)} \leq 2\varepsilon  \,.
\]
\end{theorem}
\medskip

\noindent
{\bf Proof.} Consider the mapping
\[
 S[\Phi ] (t,s)
  :=  
\Phi _0(t,s) +
G[ b^{-1} (t)\Gamma (t) F( b(t)  \Phi  )] (t,s)    \,.
\]
Due to the triangle inequality for
the $X$ norm,   for every $T_1 \in (t_0 ,T] $ we have
\begin{equation}
\label{4.2}
\parallel  S[\Phi ]   \parallel _{X(T_1)}
 \, \leq  \,
\parallel \Phi _0 \parallel _{X(T_1)}  + \parallel G[ b^{-1} (t)\Gamma (t) F( b(t)  \Phi  )]  \parallel _{X(T_1)}\,.
\end{equation}
Meanwhile,  according to Corollary~\ref{C2.2} and Proposition~\ref{P3.5},   we  derive
\begin{eqnarray}
\label{4.3}
\parallel G[ b^{-1} (t)\Gamma (t) F( b(t)  \Phi  )]  \parallel _{X(T_1)}
& \leq  &
C_N C_E \widetilde{C}_{\alpha  ,\Gamma  ,\alpha _0}(T_1)
 \|    \Phi   \|^{ \alpha +1}_{ X(T_1)} ,
\end{eqnarray}
with the function $ \widetilde{C}_{\alpha  ,\Gamma  ,\alpha _0}(T_1) $ (\ref{CbG}), the constant $C_E$ of Corollary~\ref{C2.2} and $C_N$   of Proposition~\ref{P3.5}. For the local existence ($T_1<\infty$) we choose the second case of (\ref{CbG}), while for the
global existence ($T_1=\infty$) the first case can be used as well.
\smallskip

In order to prove local solvability, 
we claim that for some $R_1> R  $ the operator $S$ maps $X({R_1 ,T_1})$,   into itself and that $S$  is a contraction provided that
$T_1-t_0   $ is sufficiently small. Indeed, for $\Phi  (x,t) \in X(R_1 ,T_1)$ inequalities  (\ref{4.2}) and (\ref{4.3}) imply
\begin{eqnarray*}
\parallel  S[\Phi ]  \parallel _{X(T_1)}
& \leq  &
R  + C_N C_E\widetilde{C}_{\alpha  ,\Gamma  ,\alpha _0}(T_1)
 \|    \Phi   \|^{ \alpha +1}_{ X(T_1)} \leq
R  + C_NC_E \widetilde{C}_{\alpha  ,\Gamma  ,\alpha _0}(T_1)  R_1 < R_1
\end{eqnarray*}
provided that $T_1-t_0   $ is sufficiently small since $\lim_{T_1 \to t_0 } \widetilde{C}_{\alpha  ,\Gamma  ,\alpha _0}(T_1)  =0$ .  To prove that  $S$ is a contraction we write
\begin{eqnarray*}
&  &
\parallel  S[\Phi ] -  S[\Psi  ] \parallel _{X(T_1)} \\
& = &
\parallel G[b^{-1} (t)\Gamma (t) \left( F( b(t)  \Phi  )- F( b(t)  \Psi   )\right) ]  \parallel _{X(T_1)}   \\
& \leq  &
C_N C_E\widetilde{C}_{\alpha  ,\Gamma  ,\alpha _0}(T_1) \parallel b^{-1} (t)\Gamma (t) \left( F( b(t)  \Phi  )- F( b(t)  \Psi   )\right) \parallel _{L^1 ([ t_0,T_1 ];L^2(S))}    \\
& \leq  &
C_N C_E\widetilde{C}_{\alpha  ,\Gamma  ,\alpha _0}(T_1) \max_{\Omega =\Phi ,\Psi } \parallel   \Omega    \parallel _{X(T)}^\alpha   \parallel  \Phi   -   \Psi  \parallel _{X(T_1)}
\end{eqnarray*}
and then choose $T_1 $ such that $\displaystyle C_NC_E\widetilde{C}_{\alpha  ,\Gamma  ,\alpha _0}(T_1) R_1^\alpha<1 $
due to $\lim_{T_1 \to t_0 } \widetilde{C}_{\alpha  ,\Gamma  ,\alpha _0}(T_1)  =0$.
\smallskip

In order to prove  global solvability ($T_1=T=\infty$),
we are going to prove that for the given $\Phi_0  (x,t) \in X(\varepsilon  ,\infty)$, the operator   $S$ maps $X({\varepsilon_1 ,\infty})$, $\varepsilon < \varepsilon_1   <\varepsilon_0 $,   into itself and that  $S$ is a contraction provided that
$\varepsilon _0  $ and $\varepsilon_1  $ are sufficiently small.
Thus,
\begin{eqnarray*}
\parallel  S[\Phi ]  \parallel _{X(\infty)}
& \leq  &
   \parallel \Phi _0\parallel _{X(\infty)}    +C_N C_E  {C}_{\alpha  ,\Gamma  ,\alpha _0}(\infty)
 \|     \Phi   \|^{ \alpha +1}_{ X(\infty)} \, ,
\end{eqnarray*}
and the operator $S$ maps  the space $X({\varepsilon_1 ,\infty})$   into itself
provided that
\begin{eqnarray*}
&  & 
\varepsilon  
+  C_N C_E {C}_{\alpha  ,\Gamma  ,\alpha _0}(\infty)  \varepsilon_1 ^{ \alpha +1} \leq   \varepsilon_1  \,.
\end{eqnarray*}
To prove that\, $S$ \,is a contraction we consider
\begin{eqnarray*}
&  &
\parallel  S[\Phi ] (x,t) -  S[\Psi  ]  \parallel _{X(\infty)} \\
& = &
\parallel G[b^{-1} (t)\Gamma (t) \left( F( b(t)  \Phi  )- F( b(t)  \Psi   )\right) ]  \parallel _{X(\infty)}   \\
& \leq  &
C_N C_E {C}_{\alpha  ,\Gamma  ,\alpha _0}(\infty)\parallel b^{-1} (t)\Gamma (t) \left( F( b(t)  \Phi  )- F( b(t)  \Psi   )\right)   \parallel _{L^1 ([ t_0,T ];L^2(S))}    \\
& \leq  &
C_N C_E {C}_{\alpha  ,\Gamma  ,\alpha _0}(\infty)\max_{\Omega =\Phi ,\Psi } \parallel   \Omega    \parallel _{X(\infty)}^\alpha   \parallel  \Phi -   \Psi  \parallel _{X(\infty)}
\end{eqnarray*}
and choose $\displaystyle C_N C_E {C}_{\alpha  ,\Gamma  ,\alpha _0}(\infty) \max_{\Omega =\Phi ,\Psi } \parallel   \Omega    \parallel _{X(\infty)}^\alpha <1 $.
Theorem is proven. \hfill $\square$
\medskip

\noindent
{\bf Proof of Theorem~\ref{T1.1}.} We have to prove that the function $\Phi _0(x,t)  $,
 generated by the Cauchy problem (\ref{IC}) for the linear equation  without source, belongs to $  X(\varepsilon ,T)$ and that it has sufficiently small norm.
Indeed, according to Corollary~\ref{C2.2}, we have the following estimate 
\begin{eqnarray*}
 &  &
    \| \partial_t \Phi _0  \| _{L^\infty([t_0 ,t]; L^2(S))} +  \| a^{-1} (\cdot )\nabla \Phi _0  \| _{L^\infty([t_0 ,t]; L^2(S))}
 + \| M (\cdot )\Phi _0  \| _{L^\infty([t_0 ,t]; L^2(S))}  \nonumber \\
&  & 
+  \|   \sqrt{ \dot{a}  a^{-3}  } \,\nabla \Phi _0   \|_{L^2([t_0 ,t]\times S)}   +  \| \sqrt{ |\dot{c} | } \, \Phi _0  \| _{L^2([t_0 ,t]\times S)}   \nonumber
 \\
 & \leq &
C_E\left(   \| \varphi _1\|_{L^2(S)}  + a^{-1}(t_0 )\|   \nabla \varphi _0   \|_{L^2(S)}    +   M (t_0 )   \|    \varphi _0  \|_{L^2(S)}   \right)
\end{eqnarray*}
for the solution of the linear problem without source, consequently,
\[
\|  \Phi _0  \|_{X(T)}
\leq
C_E\left(   \| \varphi _1 \|_{L^2(S)}  + a^{-1}(t_0 )\|   \nabla \varphi _0  \|_{L^2(S)}    +   M (t_0 )   \|    \varphi _0  \|_{L^2(S)}   \right)  .
\]
It remains to set the right hand side of the last inequality sufficiently small by the proper choice of the initial functions and to apply the statement  (ii) of
Theorem~\ref{T4.1}.  Theorem is proven.  \hfill $\square$
\medskip

\noindent
{\bf Proof of Theorem~\ref{T1.3}.} We just repeat the above argument and then apply the statement (i) of
Theorem~\ref{T4.1}.
Theorem is proven.  \hfill $\square$
\medskip

\noindent
{\bf Proof of Theorem~\ref{T1.2}.} With the potential function $V (t,x,u)  $, we define the energy
\begin{eqnarray*}
&  &
E_V(t)  \\
& :=   &
\frac{1}{2} \Big\{ \| u_{t}\|^2_{ L^2( S )}+  a^{-2}(t)\| \nabla_\sigma  u \|^2_{ L^2( S )}
 +   M^2(t)   \|      u   \|^2_{ L^2(S )} +  b^{-2}(t)  \int_{S} V (t,s,b(t)u(t,s))\,d\mu _\sigma \Big\} \,.
\end{eqnarray*}
Then we have 
\begin{eqnarray*}
\frac{d}{d t }  E_V(t )
& = &
  \frac{1}{2}  ( a^{-2}(t) )_t\| \nabla_\sigma  u \|^2_{ L^2( S )}
+  \frac{1}{2}  (M^2 (t) )_t   \|      u   \|^2_{ L^2(S )} \\
&  &
  +   \frac{1}{2}  \int_{S} b^{-2}(t)\Bigg\{ V_t(t,s,b(t)u(t,s))  -  2b^{-1}(t)\dot{b}(t)V (t,s,b(t)u(t,s))  \\
&  &
+       \dot{b} (t)  u(t,x)  V_\psi (s,b(t)u(t,s))\Bigg\} \,d\mu _\sigma  \,.
\end{eqnarray*}
Assumption (\ref{potential}) implies 
\[
 V_t(t,s,b(t)u(t,s))  -  2b^{-1}(t)\dot{b}(t)V (t,s,b(t)u(t,s))   
+       \dot{b} (t)  u(t,x)  V_\psi (s,b(t)u(t,s))\leq 0
\]
for all $t \in [t_0,\infty) $, $s \in S$, $w  \in {\mathbb R} $.
The integration gives
\begin{eqnarray*}
E_V(t)
& -   &
\int_{t_0}  ^t \Bigg\{  \frac{1}{2}  ( a^{-2}(\tau ) )_\tau \| \nabla_\sigma  u \|^2_{ L^2( S )}
+  \frac{1}{2}  (M^2 (\tau ) )_\tau    \|      u   \|^2_{ L^2(S )}  \Bigg\} \,d \tau \\
&  &
- \frac{1}{2}\int_{t_0}  ^t    \int_{S} b^{-2}(\tau )\Bigg\{  V_t(\tau ,s,b(\tau )u(\tau ,s))  -  2b^{-1}(\tau )\dot{b}(\tau )V (\tau ,s,b(\tau )u(\tau ,s))  \\
&  &
+       \dot{b} (\tau )  u(\tau ,x)  V_\psi (s,b(\tau )u(\tau ,s))\Bigg\} \,d\mu _\sigma   \,d \tau  = E_V(t_0  )\,.
\end{eqnarray*}
In particular,  due to the  assumption   (\ref{potential}) and assumptions  on  $a=a(t)$  and  $M=M(t)$,  we obtain
\[
E_V(t) -
\int_{t_0}  ^t \Bigg\{  \frac{1}{2}  ( a^{-2}(\tau ) )_\tau \| \nabla_\sigma  u \|^2_{ L^2( S )}
+  \frac{1}{2}  (M^2 (\tau ) )_\tau    \|      u   \|^2_{ L^2(S )}  \Bigg\} \,d \tau \leq  E_V(t_0  )\,.
\]
Since $ E_V(t)$ does not blow up in finite time, the local solution can be extended globally for all $t\geq t_0$. Theorem is proven. \hfill $\square$

\subsection{Examples}

To make examples more transparent, in this subsection we can restrict them to the case of the manifold $S$ with the single global chart, that is $S={\mathbb R}^n $. In that case the   Laplace-Beltrami operator $ \Delta _\sigma $ in the metric $\sigma $ can be simplified and the equation (\ref{CKGE}) can be rewritten 
as follows
\begin{eqnarray} 
\label{4.7}
 \psi _{tt}
-  \frac{1}{a^{2 }(t)\sqrt{|\sigma (x)|}}\frac{\partial }{\partial x^i}\left( \sqrt{|\sigma (x)|} \sigma ^{ik} (x)\frac{\partial \psi   }{\partial x^k} \right)+ n\frac{\dot{a}(t) }{a(t)}\psi_t  + m^2 \psi =  -  V'_\psi (t,x,\psi  ) \,,
\end{eqnarray} 
where the coefficients $\sigma ^{ik} (x) $, $i,k=1,2,\ldots,n\, $, belong to the space ${ B}^1({\mathbb R}^n) $ of the functions with the uniformly bounded derivatives. Moreover, the symmetric form $ \sigma ^{ik} (x) \xi _i\xi _k$ is  positive:  $ \sigma ^{ik} (x) \xi _i\xi _k\geq const >0$ for all  $x$, $\xi \in  {\mathbb R}^n$, $|\xi |=1 $. 

\noindent
\begin{example} \label{E4.1} \mbox{\rm Consider the equation (\ref{CKGE}) with $ V'_\psi (x,\psi  )=  -\Gamma (t)F(s,\psi )$, $
a(t)=  t^{ {\ell}/{2}} $,  $\ell >0 $:}  \end{example}  \vspace{-0.5cm} 
\begin{eqnarray*}
    \psi _{tt}
-  t^{-\ell } \Delta_\sigma      \psi  + \frac{ n \ell }{ 2 t }       \psi_t + m^2 \psi
& = &
 \Gamma (t)F(s,\psi )\,.
\end{eqnarray*}
 The number $\ell  =4/3 $\, for $n=3$ coincides with the Einstein-de~Sitter exponent.  
We make change $\psi =t^{-  n\ell/{4}}u $, then
\begin{eqnarray*}
  u_{tt} -  t^{-\ell  } \Delta_\sigma      u +  M^2(t)   u
& = &
 t^{n\ell /4} \Gamma (t)F(s,t^{-\frac{ n\ell}{4}}u )
\end{eqnarray*}
with the curved mass $M^2_{EdS}(t)$ is given by (\ref{MEdS}).
If $n\ell=4 $ and, in particular, in the case of the Einstein-de~Sitter spacetime with $n=3 $
and $l=4/3 $ the curved mass coincides with the physical mass,   and the equation is
\begin{eqnarray*}
  u_{tt}
-   t^{-{4}/{3}}   \Delta_\sigma      u    +m^2u =t \Gamma (t)F(s,t^{-1}u )\,.
\end{eqnarray*}  
The condition (\ref{Ma}) means  $  n  \ell \leq 4$. 
The condition (\ref{GammaB}) reads 
\[
|\Gamma (t)| \lesssim     t^{  -1} \quad  \mbox{\rm for all} \quad t \geq t_0\,.
\]
For (\ref{GammaU}) we have
\begin{eqnarray*}
&    &
 \int_{t_0}  ^\infty t^{\frac{n\alpha _0  }{4-n\alpha _0} }  \left| \Gamma (t)\right|^{\frac{4}{4-n\alpha _0}} \,dt   < \infty \, ,
\end{eqnarray*}
where $0 <  n\alpha _0 < 4  $. If $ \Gamma (t)=t^\gamma $, then the last integral is convergent if  $\gamma <-1 $.

\begin{example} \label{E4.2}
\mbox{\rm Consider a spacetime with the scale function}\end{example}   \vspace{-0.3cm} 
\begin{eqnarray*}
& &
a(t)= \exp(H  t^\beta ),\quad H, \beta  \in {\mathbb R}\,,
\end{eqnarray*}
and 
\begin{eqnarray*}
&  &
c(t)=  m^2  -\frac{1}{4} \beta  H n t^{\beta -2} \left(\beta  H n t^{\beta }+2\beta  -2\right), \quad   \\ 
& &
 \dot{c} (t)= 
-\frac{1}{2} (\beta -1) \beta  H n t^{\beta -3} \left(\beta +\beta  H n t^{\beta }-2\right)\,.
\end{eqnarray*}
If  $H$ is positive, then it can be regarded as  the Hubble constant.   
For the condition (\ref{a})  we have to set $  \beta H>0$. 
For the  condition (\ref{Ma})  we consider two cases:\\
 {\bf 1.  } $ \beta >0 $.  
Both conditions of (\ref{Ma}) are fulfilled only
if  $\beta = 1 $ and $0< H <\frac{2m}{n} $. That is a case of the de~Sitter metric and the large mass in the classification of \cite{Yagdjian-Galstian}. 
For $\beta = 1 $ the condition  (\ref{GammaB}) reads 
\[
|\Gamma (t)| \lesssim 1\quad  \mbox{\rm for all} \quad t \geq t_0.
\]
The condition (\ref{GammaU}) for $ \Gamma (t)=t^\gamma$  implies  $\gamma   < -1+\frac{1}{4} \alpha _0   n  $. \\
The condition (\ref{potential}) for the self-interaction $ V'_\psi (x,\psi  )=  -  \Gamma (t)\mu (x)|\psi|^\alpha \psi  $ and $\beta =1 $ is satisfied if $\dot{\Gamma } (t) \leq \frac{\alpha nH}{ 2} {\Gamma } (t)$, that is $0\leq \Gamma (t) \leq C\exp(\frac{\alpha nH}{ 2}t)  $, and $0<c_0 \leq \mu (x)\leq c_1$. We do not know if   the last condition on $\Gamma (t) $  is a necessary restriction  for the case of the energy conservative  potentials.\\ 
  {\bf 2.} $ \beta <0 $.   In this case  $H<0 $. 
 Both conditions of (\ref{Ma})  are fulfilled
for  $\beta <0$ and   $ H <0 $. \\
For the condition (\ref{GammaB}) we obtain
\[
|\Gamma (t)| \lesssim     t^{\beta -1} \quad  \mbox{\rm for all} \quad t \geq t_0
\]
while for  (\ref{GammaU}) we have
\begin{eqnarray*}
&    &
  \int_{t_0}  ^\infty    t^{(1-\beta)  \frac{n\alpha _0  }{4-n\alpha _0} }  \left| \Gamma (t)\right|^{\frac{4}{4-n\alpha _0}} \,dt < \infty  ,
\end{eqnarray*}
where $0 <  n\alpha _0 < 4  $. If $ \Gamma (t)=t^\gamma $, then for the convergence of the last integral  we need 
$
\gamma   < -1+\frac{1}{4}  \beta  n\alpha _0$. 
Hence, $\gamma  $ must  decay with the rate constant:
\[
\gamma \leq   \beta -1 \quad or \quad \gamma < \frac{n\alpha _0}{4}  \beta  -1\,.
\]
We do not know if   the last condition is a necessary restriction.

\begin{example}  
\mbox{\rm Consider the following example of the scale function} \end{example}   \vspace{-0.3cm} 
\begin{eqnarray*}
& &
a(t)=  t^\frac{\ell}{2} \exp(H t^{\beta } ),  \qquad \ell, H, \beta  \in {\mathbb R},\\
\end{eqnarray*}
and
\begin{eqnarray*}
&  &
\dot{a}(t)= t^{\frac{\ell}{2} -1}  \left(\beta  H t^{\beta }+\frac{\ell}{2} \right)\exp(H t^{\beta } ),    \\ 
 &  &
c(t)=  m^2-\frac{n \left(\beta ^2 H^2 n t^{2 \beta }+2 \beta  H t^{\beta } (\beta +\frac{\ell}{2}  n-1)+\frac{\ell}{2}  (\frac{\ell}{2}  n-2)\right)}{4 t^2}>0\,\,  \mbox{\rm for  } \,\, t\geq t_0\,,\\
&  &
 \dot{c} (t)= -\frac{n \left((\beta -1) \beta ^2 H^2 n t^{2 \beta }+(\beta -2) \beta  H t^{\beta } (\beta +\frac{\ell}{2}  n-1)+\lambda  (2-\frac{\ell}{2}  n)\right)}{2 t^3}\leq 0  \,\,  \mbox{\rm for  } \,\, t\geq t_0.
 \end{eqnarray*}
Here  $t_0$ is sufficiently large number. The case of $\ell =0 $  coincides with  Example~\ref{E4.2}, while the case of  $\beta  =0$ or $H=0 $ coincides with  Example~\ref{E4.1}.

\noindent For the conditions (\ref{a})-- (\ref{GammaU}) we consider two separate cases.\\
\noindent
{\bf 1.  }  $\beta >0 $.  In this case the condition  $\dot{a}(t)>0 $ of (\ref{a})  implies $H>0 $.\\
The condition $c(t)>0 $ of (\ref{Ma}) is satisfied for $ \beta <1 $ or $\beta =1  \,\, \&\,\, m^2  -\frac{1}{4}    (H n)^2  >0$, 
while the condition $\dot{c} (t)\leq 0 $ implies $ \beta >1 $ or $\beta =1 $ and $H=0$.
Thus, (\ref{Ma}) is satisfied only for $\beta =1  $ and $H=0$, which coincides with  Example~\ref{E4.1}.
\smallskip

\noindent
{\bf 2.  }  $\beta <0 $. In this case the assumption $\dot{a}(t)>0 $ of (\ref{a}) implies $\ell >0 $.
Both conditions of (\ref{Ma}) are satisfied  for $0\leq \ell \leq \frac{4}{n} $. Now, if $\ell >0 $ for the condition (\ref{GammaB}) 
we obtain 
\[
|\Gamma (t)| \lesssim  t^{-1} \quad  \mbox{\rm for all} \quad t \geq t_0 \,,
\]
and for  (\ref{GammaU}) we have
\[
  \int_{t_0} ^\infty  \left( \frac{t}{\beta  H t^{\beta }+\frac{\ell}{2}}\right)^\frac{n\alpha _0  }{4-n\alpha _0}   \left| \Gamma (t)\right|^{\frac{4}{4-n\alpha _0}} \,dt < \infty  ,
\]
where $0 <  n\alpha _0 < 4  $. If $ |\Gamma (t)| =t^\gamma $, then for the convergence of the last integral  we need 
$\gamma   < -1$. Hence, the condition for the decay rate constant $\gamma $  coincides  with  the  condition   in  Example \ref{E4.1}\,.
\medskip

Finally, we note here that the results of this paper are applicable to the equation (\ref{4.7}) with $x$-dependent coefficient, 
while the results of  \cite{Yag_2005},\cite{JMAA_2012},\cite{Helsinki_2013},\cite{Nakamura} are restricted to the equation with $x$-independent coefficients.

\smallskip

\section*{Acknowledgments}
The authors are indebted to the anonymous referee for for pointing out the  discrepancy in the initial 
version of Lemma~\ref{L3.1} and for the numerous remarks, comments and suggestions which improved  the readability of the text.
The referee  also  drew  our attention to the articles \cite{Ringstrom,Rodnianski-Speck}.


\begin{thebibliography}{99}


\bibitem{Baskin}
D.~Baskin,~A parametrix for the fundamental solution of the Klein-Gordon equation on asymptotically
de Sitter spaces,  Journal of Functional Analysis ~{259}~(7)~(2010)~1673--1719.

\bibitem{BaskinSE}
D.~Baskin,~Strichartz estimates on asymptotically de Sitter spaces, Ann. Henri Poincar\'e 14 (2) (2013)~221-252.


  
\bibitem{Choquet-Bruhat_book}
Y.~Choquet-Bruhat,   General relativity and the Einstein equations. Oxford Mathematical Monographs. Oxford University Press, Oxford, 2009.

\bibitem{Choquet-Bruhat-ND_2000}
Y.~Choquet-Bruhat, Global wave maps on Robertson-Walker spacetimes. Modern group analysis. Nonlinear Dynam. 22 (1) (2000)~39--47.

 

\bibitem{Hawking}
S.W.~Hawking,~ G.F.R.~Ellis,~The large scale structure of space-time. Cambridge Monographs on Mathematical Physics,
No. 1. Cambridge University Press, London-New York, 1973.


\bibitem{Higuchi}
A.~Higuchi, Forbidden mass range for spin-$2$ field theory in de~Sitter spacetime, Nuclear Phys. B 282 (2) (1987)  397--436.


\bibitem{Hintz-Vasy}
P.~Hintz, A.~Vasy, Semilinear wave equations on asymptotically de Sitter, Kerr-de Sitter and Minkowski spacetimes, 	arXiv:1306.4705.


\bibitem{Nakamura}
M.~Nakamura,~The Cauchy problem for semi-linear Klein–Gordone quations in de Sitter spacetime. J. Math. Anal. Appl. { 410}, Issue 1,   (2014)   445-454.



\bibitem{Polya-Szego}
G.~P\'olya, G.~Szeg{\"o},   Problems and theorems in analysis. I. Series, integral calculus, theory of functions.  Grundlehren der Mathematischen Wissenschaften, 193. Springer-Verlag, Berlin-New York, 1978.

\bibitem{Ringstrom}
H.~Ringstr\"om,   Future stability of the Einstein-non-linear scalar field system. Invent. Math.  173,  no. 1  (2008)  123--208.

\bibitem{Rodnianski-Speck}
I.~Rodnianski,   J.~Speck,   The nonlinear future stability of the FLRW family of solutions to the irrotational Euler-Einstein system with a positive cosmological constant. J. Eur. Math. Soc. (JEMS)  15,  no. 6 (2013)  2369--2462.


\bibitem{Shatah}
J.~Shatah,~ M.~Struwe, 
{ Geometric wave equations.}
Courant Lecture Notes in Mathematics, 2. New York University,
Courant Institute of Mathematical Sciences, New York; American Mathematical Society, Providence, RI, 1998.




\bibitem{Tartar}
L.~Tartar,~ An Introduction to Sobolev Spaces and Interpolation Spaces, Springer, 2007.

\bibitem{Vasy_2010}
A.~Vasy,~The wave equation on asymptotically de Sitter-like spaces.~{Adv. Math.}  {223} (1) (2010)~49--97.



\bibitem{Yag_2005}
K.~Yagdjian,~
Global existence in the Cauchy problem for nonlinear wave equations with variable speed of propagation.   New trends in the theory of hyperbolic equations, 301-385,
Oper. Theory Adv. Appl., 159, Birkh\"auser, Basel, 2005.





\bibitem{Yagdjian-Galstian}
K.~Yagdjian,~A.~Galstian,  Fundamental solutions for the Klein-Gordon equation in de Sitter spacetime.
Comm. Math. Phys. {  285} (2009)~ 293--344. 

\bibitem{JMAA_2012}
K.~Yagdjian,~Global existence of the scalar field in de Sitter spacetime. J. Math. Anal. Appl. 396~(1) (2012)  323--344.



\bibitem{Helsinki_2013}
K.~Yagdjian,~
Semilinear Hyperbolic Equations in Curved Spacetime. Fourier Analysis,
Pseudo-differential Operators, Time-Frequency Analysis and Partial Differential Equations.
Series: Trends in Mathematics. Birkh\"auser Mathematics,   391--415, 2014.


\end{thebibliography}
\end{document}